%
%
%
\magnification=1200
\documentstyle{amsppt}
\def\supp{\text{supp}}
\def\M{\text{M}}
\def\C{\text{C}}
\def\MC{\text{M}_{\Bbb C}}
\def\L{{\Cal L}}
\def\uto{\rightrightarrows}
\topmatter
\title
Relative asymptotics for polynomials orthogonal with respect to
 a discrete Sobolev inner product
\endtitle
\author
Guillermo L\'opez, Francisco Marcell\'an and Walter Van Assche
\endauthor
\affil
Universidad de la Habana \\ Universidad Carlos III de Madrid \\
Katholieke Universiteit Leuven
\endaffil
\address
(G.L.) Facultad de Matem\'atica y Cibern\'etica,
Universidad de la Habana,
Habana~4 (Cuba)
\endaddress
\address
(F.M.) Departamento de Ingenier\'\i a,
Escuela Polit\'ecnica Superior,
Universidad Carlos III de Madrid,
E-28913 Legan\'es (SPAIN)
\endaddress
\email pacomarc\@ing.uc3m.es  \endemail
\address
(W.V.A.) Katholieke Universiteit Leuven,
Department of Mathematics,
Ce\-les\-tijn\-en\-laan 200\,B,
B-3001 Heverlee (BELGIUM)
\endaddress
\email fgaee03\@cc1.KULeuven.ac.be \endemail
\thanks
This research was carried
out will the first author was visiting ETSII, Universidad Polit\'ecnica
de Madrid, under sabbatical grant 75/91 of Secretar\'\i a de Estado de
Universidades e Investigaci\'on of the Ministerio Espa\~nol de
Educaci\'on y Ciencia.
The second author was supported by Comisi\'on Interministerial de Ciencia
y Tecnolog\'\i a (CICYT) PB-89/0181/C02/01.
The third author is a Research Associate of the Belgian National Fund
for Scientific Research.
\endthanks
\keywords
Orthogonal polynomials, Sobolev inner product, asymptotic behavior
\endkeywords
\subjclass
42C05
\endsubjclass
\rightheadtext{Relative asymptotics for Sobolev inner products}
\leftheadtext{L\'opez, Marcell\'an and Van Assche}
\abstract
We investigate the asymptotic properties of orthogonal
polynomials for a class of inner products including the discrete Sobolev inner
products $\langle h,g \rangle = \int hg\, d\mu
+ \sum_{j=1}^m \sum_{i=0}^{N_j} M_{j,i} h^{(i)}(c_j) g^{(i)}(c_j)$,
where $\mu$ is a certain type of
complex measure on the real line, and $c_j$ are complex numbers in the
complement of $\supp(\mu)$. The Sobolev orthogonal polynomials are compared
with the orthogonal polynomials corresponding to the measure $\mu$.
\endabstract

\endtopmatter

\document

\head 1. Introduction \endhead

In the last few years orthogonal polynomials on Sobolev spaces have attracted
considerable attention. This interest is justified for several reasons:

\roster
\item The comparison with the standard theory of orthogonal polynomials;
\item The spectral theory of ordinary differential equations;
\item The analysis of spectral methods in the numerical treatment of partial
differential equations;
\item The search of algorithms for computing Fourier Sobolev series in terms of
Sobolev orthogonal polynomials for norms involving the function and some of its
derivatives.
\endroster

Despite the effort, most of the progress attained is for special types of inner
products. A survey on the subject is provided in \cite{2}. The results known
are mostly connected with the formal theory: recurrence relations, location of
zeros, differential formulas and so on. Little is known concerning asymptotic
properties. Possibly, the only results known in this direction are contained in
\cite{15} for a special type of so called discrete Sobolev inner product.

Let $\mu$ be a  finite  positive  Borel measure  whose  support, $\supp(\mu) =
S_\mu$, contains an infinite set of points, $S_\mu  \subset  {\Bbb R}$.
A discrete Sobolev inner product is given by the expression
  $$   \langle h,g \rangle = \int h g \, d\mu +
    \sum_{j=1}^m \sum_{i=0}^{N_j} M_{j,i} h^{(i)}(c_j) g^{(i)}(c_j) ,
                        \tag 1.1  $$
where $c_j \in  {\Bbb R}$,  $M_{j,i} \geq  0$,  $m, N_j > 0$
(see \cite{5}).

By $S_n$, $n \in  {\Bbb Z}_+=\{1,2,3,\ldots\}$,
we denote the monic polynomial of least degree,
not identically zero, such that
  $$  \langle p, S_n \rangle = 0, \qquad  p \in {\Bbb P}_{n-1},  \tag 1.2 $$
where ${\Bbb P}_{n-1}$ is the linear space of all polynomials of degree
less than or equal to
$n-1$. It is easy to see that for every $n$ the degree of $S_n$ is $n$.
Furthermore let $L_n$
be the monic orthogonal polynomials of degree $n$ with respect to $\mu$.
When \thetag{1.1} reduces to
  $$  \langle h,g \rangle = \int hg\, d\mu + M_1 h'(c) g'(c), \qquad M_1 > 0,
     \tag 1.3 $$
where $\mu$ belongs to the  class $\M(0,1)$
(for the definition see \cite{19}, but also \thetag{2.1} below) and $c \in
{\Bbb R} \setminus S_\mu$, then, as proved in \cite{15},
  $$  \lim_{n \to \infty} \frac{S_n(z)}{L_n(z)} =
     \frac{(\varphi(z)-\varphi(c))^2}{2\varphi(z) (z-c)},
     \tag 1.4 $$
uniformly on every compact subset of $\overline{\Bbb C} \setminus S_\mu$,
where
  $$ \varphi(z) = z + \sqrt{z^2-1}, \qquad
    (\sqrt{z^2-1} > 0 \text{ for } z > 1). $$
In \cite{15} an equivalent expression is given for the right-hand side of
\thetag{1.4}.
Throughout the paper we will use the notation $\varphi_n \uto \varphi,\
K \subset D$ when the sequence of functions $\varphi_n$ converges to
$\varphi$ uniformly on every compact subset $K$ of the region $D$. In this
way \thetag{1.4} becomes
  $$   \frac{S_n(z)}{L_n(z)} \uto
     \frac{(\varphi(z)-\varphi(c))^2}{2\varphi(z) (z-c)},
     \qquad K \subset \overline{\Bbb C} \setminus S_\mu .  $$

The denominators $Q_n$ of the main diagonal sequence for Pad\'e approximants of
Stieltjes type meromorphic functions
  $$  f(z) =  \int \frac{d\mu(x)}{z-x} + \sum_{j=1}^m \sum_{i=0}^{N_j}
    A_{j,i} \frac{i!}{(z-c_j)^{i+1}}, \qquad A_{j,N_j} \neq 0,
    \tag 1.5 $$
satisfy orthogonality relations similar to \thetag{1.2}. In fact
  $$  0 = \int p Q_n \, d\mu +  \sum_{j=1}^m \sum_{i=0}^{N_j}
     A_{j,i} \left( p(z) Q_n(z) \right)^{(i)}_{z=c_j},
     \qquad p \in {\Bbb P}_{n-1}. \tag 1.6 $$
Here $A_{j,i}$ may be complex numbers, and of particular interest is the
case
when the points $c_j$ belong to the complement of the convex hull of $S_\mu$.
Of course, the inner product associated with \thetag{1.6} is not positive
definite.
Therefore, a priori, it is not possible to guarantee that the
degree of $Q_n$ is $n$.
Nevertheless, in \cite{8} A.A. Gonchar proves the following result.
Assume that $S_\mu \subset [-1,1]$, and $\mu$ is such that
$$   \frac{L_{n+1}(z)}{L_n(z)} \uto
   \frac{\varphi(z)}2, \qquad K \subset  {\Bbb C} \setminus [-1,1],
     \tag 1.7 $$
then
 $$   \frac{Q_n(z)}{L_n(z)} \uto
  \prod_{j=1}^m \left( \frac{(\varphi(z)-\varphi(c_j))^2}{2\varphi(z) (z-c_j)}
   \right)^{N_j+1}, \qquad K \subset \overline{\Bbb C} \setminus [-1,1].
   \tag 1.8 $$

A well known result (see \cite{19}) asserts that \thetag{1.7} with $[-1,1]$
replaced by
$S_\mu$ is equivalent to $\mu \in  \M(0,1)$. In this case, $S_\mu$ consists
of $[-1,1]$ and at most a denumerable set of mass points contained in
${\Bbb R} \setminus [-1,1]$ which, when infinitely many, may
accumulate only at the extreme points $\pm  1$.
The best known sufficient condition for $\mu  \in  \M(0,1)$ is due to E.A.
Rakhmanov published in two  papers \cite{20} \cite{21}. This condition is
$\mu '>0\ a.e$. on $[-1,1]$ and $S_\mu  = [-1,1]$. Simplified proofs may be
found in \cite{16} and \cite{22}.

When there is one simple pole in \thetag{1.5} $(m=1, N_1 = 0)$ and $S_\mu  =
[-1,1]$, then
 the right-hand sides of the asymptotic formulas \thetag{1.4} and \thetag{1.8}
coincide
although the corresponding inner products  differ
in the atomic aggregate. In the second case, the inner product is given by
  $$  \langle h,g \rangle = \int h g \, d\mu + M_0 h(c) g(c). \tag 1.9 $$
The agreement is deeper than formulas \thetag{1.4} and \thetag{1.8} suggest.
In \cite{19, Lemma 16, \S 7}, a detailed study is carried out comparing the
asymptotic behavior of the orthonormal polynomials
(and their leading coefficients) with respect to \thetag{1.9} (with $M_0>0$)
with those of $\mu$, for
all $c \in  {\Bbb R}$. As proved in \cite{15} these relations are the same
if we replace \thetag{1.9} by \thetag{1.3}.

We will consider inner products for which \thetag{1.1} and
\thetag{1.6} are particular cases (see \thetag{4.1} and \thetag{4.3}). The
corresponding
orthogonal polynomials satisfy relative asymptotic formulas of the type
given by \thetag{1.8}. In
our approach, which follows that of A.A. Gonchar in \cite{8}, formula
\thetag{1.7} applied
at the points $c_j$ plays a key role. Thus, in the following, we assume that
$c_j \in  {\Bbb C} \setminus S_\mu$ and $M_{j,i} \in  {\Bbb C}$ and
thus allow complex values for these quantities.
We will also allow $\mu $ to be a certain type of complex measure which extends
the class $\M(0,1)$. These measures are introduced in Section 3 where the
asymptotic behavior of the corresponding orthogonal polynomials is studied.
In Section 4 the main result of type \thetag{1.8} is obtained (Theorem 4).
 From it  we can obtain the following result.
Let $S_n$ be the $n$th monic orthogonal polynomial with respect to \thetag{1.1},
where
$\mu  \in  \M(0,1)$, $M_{j,i} \in  {\Bbb C}$, $c_j \in  {\Bbb C}
\setminus S_\mu$. For each fixed $j=1,\ldots, m$, we denote by $N^*_j$ the
number of coefficients $M_{j,i}$
different from zero $(i = 0, \ldots, N_j)$. If $L_n$ is the $n$th monic
orthogonal polynomial with respect to $\mu$, then
  $$  \frac{S_n(z)}{L_n(z)} \uto \prod_{j=1}^m
    \left( \frac{(\varphi(z)-\varphi(c_j))^2}
         {2 \varphi(z) (z-c_j)} \right)^{N_j^*} ,
   \qquad K \subset \overline{\Bbb C} \setminus S_\mu .     \tag 1.10 $$
Using Rouch\'e's theorem we see that for all sufficiently large $n$,
the degree of $S_n$ is $n$. Moreover, \thetag{1.1} implies that any neighborhood
of $S_\mu$
compactly contained in ${\Bbb C} \setminus  \{c_1, \ldots, c_m\}$ has, for
all sufficiently large $n$, exactly $n - N$ $(N = N^*_1 +...+ N^*_m)$ zeros of
$S_n$ (for short, we say that $n-N$ zeros of $S_n$ concentrate or accumulate on
$S_\mu$), while each sufficiently small neighborhood of $c_j$, $j =
1,\ldots,m$
contains, for all sufficiently large $n$, exactly $N^*_j$ zeros of $S_n$
(again, for short, we say that $c_j$ attracts $N^*_j$ zeros of $S_n$).

We shall limit ourselves to measures with compact support, but the methods
employed in \cite{10}, \cite{11} allow to extend the results in Section 4 to
measures with unbounded support.

Section 2 has an auxiliary character. Most of the formulas given there are
known
to the specialist, though precise reference, in the extension stated, may be
hard to find in the literature. Therefore, for completeness, we include this
short section with full proofs except when exact reference is available.

\head 2. Auxiliary results \endhead

Let $\mu$ be a finite positive measure on ${\Bbb R}$. Assume that $\mu
\in \M(0,1)$, that is, if $l_n(z) = \tau_n z^n + \cdots$, $\tau_n>0$, is the
$n$th orthonormal polynomial with respect to $\mu$ with recurrence relation
$$  z l_n(z) = a_{n+1} l_{n+1}(z) + b_n l_n(z) + a_n l_{n-1}(z), $$
then
$$   \lim_{n \to \infty} a_n = \frac12, \quad
   \lim_{n \to \infty} b_n = 0. \tag 2.1 $$
As pointed out, $\mu  \in  \M(0,1)$ is equivalent to
$$   \frac{l_{n+1}(z)}{l_n(z)} \uto \varphi(z) =
   z + \sqrt{z^2-1}, \qquad K \subset {\Bbb C} \setminus S_\mu ,
   \tag 2.2 $$
and $S_\mu$ has the structure described in Section~1.
If $L_n$ is the $n$th monic orthogonal polynomial, one also has
$$   \frac{L_{n+1}(z)}{L_n(z)} \uto \frac{\varphi(z)}2 ,
        \qquad K \subset {\Bbb C} \setminus S_\mu ,
     \tag 2.3 $$
and
$$  \lim_{n \to \infty}  \frac{\tau_{n+1}}{\tau_n} = 2 . \tag 2.4 $$
From \thetag{2.2} or from \thetag{2.1} it follows that for all fixed
$\nu \in {\Bbb Z}_+$
$$  \lim_{n \to \infty} \frac{l_{n+1}^{(\nu)}(z)}{l_n^{(\nu)}(z)}
  = 2  \lim_{n \to \infty} \frac{L_{n+1}^{(\nu)}(z)}{L_n^{(\nu)}(z)}
  = \varphi(z), \tag 2.5 $$
uniformly on compact sets of ${\Bbb C} \setminus S_\mu$.
See Lemma 4 in \cite{6} where \thetag{2.5}
is obtained from \thetag{2.2} by induction on $\nu$.
Actually, in \cite{6} it is assumed that $S_\mu \subset [-1,1]$ but the
more complicated structure of $S_\mu$ for $\mu \in \M(0,1)$ offers no real
difficulty since we know that between two consecutive points in
$S_\mu \setminus (-1,1)$, for all sufficiently large $n$, there lies
exactly
one zero of $L_n$ which  is attracted  by the point further away from $[-1,1]$
(the same is true for $l_n^{(\nu )}$).

In \cite{19} (see Lemma 1.5 of Section 4.1), P. Nevai proves that for $\mu
\in \M(0,1)$
$$   \frac1n \frac{l_n'(z)}{l_n(z)} \uto
  \frac{1}{\sqrt{z^2-1}} , \qquad K \subset \overline{\Bbb C} \setminus S_\mu.
   \tag 2.6 $$
In fact, \thetag{2.6} follows from the weaker assumption on $\mu$ that it be
such that
$$  | l_n(z)|^{1/n} \uto |\varphi(z)|, \qquad K \subset {\Bbb C} \setminus
     S_\mu. \tag 2.7 $$
Obviously \thetag{2.2} and thus \thetag{2.1} implies \thetag{2.7}. Also,
\thetag{2.7} is satisfied
if $S_\mu  = [-1,1]$ and $\C (\{\mu'>0\}) = 1/2$ where $\C(\cdot)$ is the
logarithmic capacity of the indicated set. This statement is contained in a
little known paper of P.P. Korovkin \cite{9}. Actually, \thetag{2.7} follows
from a still
weaker condition known as Ullman's criterion given in \cite{26} (see also
Chapter
4 in \cite{25}) in terms of minimal carrier capacity or an equivalent criterion
of Widom \cite{28}. Orthonormal polynomials for which \thetag{2.7} holds are
said to have regular exterior asymptotic behaviour on $[-1,1]$.

We have:
\proclaim{Lemma 1}
Let $\mu$ be such that \thetag{2.7} holds. Then, for all $\nu  \in
{\Bbb Z}_+$
$$    |l_n^{(\nu)}(z)|^{1/n} \uto |\varphi(z)|,
    \qquad K \subset {\Bbb C} \setminus S_\mu,  \tag 2.8 $$
and
$$   \frac1n \frac{l_n^{(\nu+1)}(z)}{l_n^{(\nu)}(z)}  \uto
   \frac{1}{\sqrt{z^2-1}}, \qquad K \subset \overline{\Bbb C} \setminus S_\mu.
   \tag 2.9 $$
\endproclaim

\demo{Proof}
From \thetag{2.7} it follows that the zeros of $l_n$ concentrate
on $S_\mu$. From this and the fact that the zeros of $l_n$ and $l_n'$
interlace it
readily follows that the family of functions
$\left\{ n^{-1} l_n'/l_n \right\}_{n
\geq n_0}$ is normal in $\overline{\Bbb C} \setminus S_\mu$.
Therefore, to prove \thetag{2.9} for
$\nu  = 0$ it suffices
to show that any convergent subsequence in $\overline{\Bbb C} \setminus S_\mu$
has  $1/\sqrt{z^2-1}$  as
its pointwise limit, say on ${\Bbb R}_+ \setminus S_\mu$.

Assume that
$$  \lim_{n \to \infty,\ n \in \Lambda} \frac1n \frac{l_n'(z)}{l_n(z)}
    = \psi_\Lambda(z), \qquad \Lambda \subset {\Bbb Z}_+ , \tag 2.10 $$
for $z \in K \subset \overline{\Bbb C} \setminus S_\mu$.
Let $z = x > x_0 = \max \{ t : t \in S_\mu \}$, then
$l_n(x) > 0$ and $\varphi(x) > 0$. Hence
from \thetag{2.7} (\thetag{2.8} for $\nu  = 0$) one has
$$  \lim_{n \to \infty} l_n^{1/n}(x)  = \varphi(x) $$
for $x \in K \subset (x_0,\infty)$,
or equivalently
$$ \lim_{n \to \infty}  \frac1n \log l_n(x) =  \log \varphi(x), \tag 2.11 $$
for $x \in K \subset (x_0,\infty)$.
The derivative on the left hand member of \thetag{2.11} is
equal to $n^{-1} l_n'/l_n$.
Thus, from \thetag{2.10} and \thetag{2.11} follows that
$$  \psi_\Lambda(x) = (\log \psi(x) )' = \frac{1}{\sqrt{x^2-1}}, $$
and we have \thetag{2.9} for $\nu  = 0$.

Assume that \thetag{2.8}--\thetag{2.9} hold for $\nu  = k$, then we show
that they are also true for $\nu  = k + 1$.
From \thetag{2.9} for $\nu  = k$ follows, taking $n$th roots, that
$$  \lim_{n \to \infty} \left| \frac1n \frac{l_n^{(k+1)}(x)}{l_n^{(k)}(x)}
   \right|^{1/n} =
    \lim_{n \to \infty} \left|  \frac{l_n^{(k+1)}(x)}{l_n^{(k)}(x)}
   \right|^{1/n} = 1, $$
uniformly for $x$ on compact subsets of ${\Bbb C} \setminus S_\mu$.
This, together with \thetag{2.8} for $\nu  = k$, gives \thetag{2.8} for
$\nu= k+1$.

Now, to prove  \thetag{2.9}, note that the family of functions
$\left\{ n^{-1} l_n^{(k+2)}/l_n^{(k+1)} \right\}_{n \geq n_0}$
is normal in $\overline{\Bbb C} \setminus S_\mu$.
 With this and \thetag{2.8} for $\nu = k+1$ one obtains \thetag{2.9} for $\nu
=k+1$ following the same arguments as above for $\nu  = 0$. \qed
\enddemo

\proclaim{Lemma 2}
Let $\mu  \in  \M(0,1)$, then for all $\nu  \in  {\Bbb Z}$
$$   \int \frac{l_{n+\nu}(x)l_n(x)}{z-x} \, d\mu(x) \uto
       \frac{1}{\varphi^{|\nu|}(z)\sqrt{z^2-1}},
  \qquad K \subset \overline{\Bbb C} \setminus S_\mu .     \tag 2.12 $$
\endproclaim

\demo{Proof}
Obviously, it is sufficient to consider $\nu \in {\Bbb Z}_+$. By the
Cauchy-Schwarz inequality one has for $z \in K \subset
\overline{\Bbb C} \setminus S_\mu$
$$  \left| \int \frac{l_{n+\nu}(x)l_n(x)}{z-x} \, d\mu(x) \right|
   \leq \frac{1}{d(K,S_\mu)} < \infty , $$
where $d(K,S_\mu)$ is the Euclidean distance between the two sets.
Therefore, for each fixed $\nu \in {\Bbb Z}_+$, the family of functions on the
left hand side of \thetag{2.12} is normal and uniform convergence follows from
pointwise convergence.
The pointwise limit follows from a result by Nevai \cite{18, Theorem 13 on
p.~45}. From it we have that
$$  \lim_{n \to \infty} \int  \frac{l_{n+\nu}(x)l_n(x)}{z-x} \, d\mu(x)
  = \frac1\pi \int_{-1}^1 \frac{t_\nu(x)}{z-x} \frac{dx}{\sqrt{1-x^2}} , $$
where $t_\nu$ is the $\nu$th Chebyshev orthonormal polynomial of the
first kind.
Thus \thetag{2.12} holds if we show that
$$  \frac1\pi \int_{-1}^1 \frac{t_\nu(x)}{z-x} \frac{dx}{\sqrt{1-x^2}}
    = \frac{1}{\varphi^\nu(z) \sqrt{z^2-1}} . \tag 2.13 $$
Recall that $t_0(x)=1$, $t_1(x) = x$ and for $\nu \geq 1$
$$   2x t_\nu(x) = t_{\nu+1}(x) + t_{\nu-1}(x), $$
or what is the same
$$   t_{\nu+1}(x) = 2x t_\nu(x) - t_{\nu-1}(x).  \tag 2.14  $$
Again we proceed by induction. For $\nu =0$ the formula is obtained from
Cauchy's integral formula and the residue Theorem.
For $\nu=1$ one has
$$ \align
\frac1\pi \int_{-1}^1 \frac{t_1(x)}{z-x} \frac{dx}{\sqrt{1-x^2}}
  &= \frac{z}\pi \int_{-1}^1 \frac{1}{z-x} \frac{dx}{\sqrt{1-x^2}}
  - \frac1\pi \int_{-1}^1  \frac{dx}{\sqrt{1-x^2}}  \\
  &= \frac{z}{\sqrt{z^2-1}} - 1 \\
  &= \frac{1}{\varphi(z)\sqrt{z^2-1}} .
  \endalign $$
Now assume that \thetag{2.13} holds for $\nu  = 0,1,\ldots,k$, $k\geq 1$,
then we prove that
it also holds for $\nu  = k+1$. In fact, from \thetag{2.14} and the
induction hypothesis one has
$$ \align
  \frac1\pi \int_{-1}^1 \frac{t_{k+1}(x)}{z-x} \frac{dx}{\sqrt{1-x^2}}
  &= \frac1\pi \int_{-1}^1 \frac{2x t_k(x)}{z-x} \frac{dx}{\sqrt{1-x^2}}
 - \frac1\pi \int_{-1}^1 \frac{t_{k-1}(x)}{z-x} \frac{dx}{\sqrt{1-x^2}}  \\
 &= \frac{2z}\pi \int_{-1}^1 \frac{t_k(x)}{z-x} \frac{dx}{\sqrt{1-x^2}}
  - \frac2\pi \int_{-1}^1 \frac{t_\nu(x) \, dx}{\sqrt{1-x^2}}   \\
  &\quad
  - \frac1\pi \int_{-1}^1 \frac{t_{k-1}(x)}{z-x} \frac{dx}{\sqrt{1-x^2}} \\
 &= \frac{1}{\varphi^{k-1}(x)\sqrt{z^2-1}}
  \left( \frac{2z}{\varphi(z)} - 1 \right) \\
   &= \frac{1}{\varphi^{k+1}(x)\sqrt{z^2-1}}  ,
   \endalign $$
which we wanted to prove. \qed
\enddemo

\head 3. Relative asymptotics for certain complex measures \endhead

As in Section 2, $\mu$ is a finite positive Borel measure, $\mu  \in
\M(0,1)$, $S_\mu  = \supp(\mu)$. Let $r = S/T$, after canceling
out common factors, where
$$  S(z) = \prod_{i=0}^h (z-c_i)^{A_i}, \quad
    T(z) = \prod_{j=1}^\ell (z-d_j)^{B_j}, \qquad c_i, d_j \in
    {\Bbb C} \setminus S_\mu, A_i, B_j \in {\Bbb N}. $$
Set
$$  A = A_1 + \cdots + A_h , \quad B = B_1 + \cdots + B_\ell. $$
Assume that $Q_n$ is the monic polynomial of least degree, not identically
equal to zero, such that
$$ 0 = \int p(x) Q_n(x) r(x) \, d\mu(x), \qquad   p \in  {\Bbb P}_{n-1},
         \tag 3.1 $$
and $L_n$ the $n$th monic orthogonal polynomial with respect to $\mu$. In the
past few years, the relative asymptotic behavior of the polynomials arising
from a general real modification of the measure $\mu$ has been considered, see
e.g., \cite{17}, \cite{18}, \cite{19} and \cite{22}. We are interested in the
asymptotic behavior of $\{ Q_n/L_n \}$ $n \in {\Bbb Z}_+$. The
corresponding result (formula (3.2) below with $\nu  = 0$) was stated as
Theorem 10 in \cite{12} but the proof was omitted being similar to that of
Theorem 6 of the same paper. Since we also need to derive other asymptotic
formulas and
there are minor changes in the assumptions, we include here the proof.

\proclaim{Theorem 1}
 Let $\mu  \in  \M(0,1)$. Then for all sufficiently large $n$ the degree
of $Q_n$ is $n$ and for all fixed $\nu  \in  {\Bbb Z}_+$
$$  \frac{Q_n^{(\nu)}(z)}{L_n^{(\nu)}(z)} \uto
   \left( \frac12 \right)^A \prod_{j=1}^\ell \left(
     1 - \frac{1}{\varphi(z) \varphi(d_j)} \right)^{B_j}
     \prod_{i=1}^h \left( \frac{\varphi(z)-\varphi(c_i)}{z-c_i} \right)^{A_i},
     \qquad K \subset \overline{\Bbb C} \setminus S_\mu, \tag 3.2 $$
and
$$ \frac{Q_{n+1}^{(\nu)}(z)}{Q_n^{(\nu)}(z)} \uto \frac{\varphi(z)}{2},
   \qquad K \subset \Bbb C \setminus S_\mu, \tag 3.3  $$
and
$$   \frac1n
\frac{Q_n^{(\nu+1)}(z)}{Q_n^{(\nu)}(z)} \uto
\frac{1}{\sqrt{z^2-1}},
   \qquad K \subset \overline{\Bbb C} \setminus S_\mu .  \tag 3.4  $$
\endproclaim

\demo{Proof}
 First we concentrate on \thetag{3.2} for $\nu  = 0$. Taking $p(x) =
T(x)q(x)$, where $q$ has degree $\leq n - B - 1$, we have
$$  0 = \int q(x) Q_n(x) S(x) \, d\mu(x). $$
Therefore
$$  R_n(z) := S(z) Q_n(z) = \sum_{k=0}^{A+B} \lambda_{n,k} L_{n+A-k}(z),
    \tag 3.5 $$
where either $\lambda_{n,0} = 1$ or  $Q_n$ has degree $< n$. Since $R_n$ is
monic, the first coefficient $\lambda_{n,k}$
different from zero appearing in \thetag{3.5} must be one.
Dividing this relation by $L_{n-B}$ we get
$$   \Omega_n(z) := \frac{R_n(z)}{L_{n-B}(z)} =
    \sum_{k=0}^{A+B} \lambda_{n,k} \frac{L_{n+A-k}(z)}{L_{n-B}(z)} .  $$
Set $\lambda_n^* = \left( \sum_{k=0}^{A+B} |\lambda_{n,k}| \right)^{-1}
< \infty$ and introduce the polynomials
$$  p_n(z) = \sum_{k=0}^{A+B} \lambda_{n,k} z^{A+B-k}, \quad
    p_n^*(z) = \lambda_n^* p_n(z) . $$
We will prove that
$$  p_n(z) \uto p_0(z) = \prod_{i=1}^h \left( z - \frac{\varphi(c_i)}{2}
  \right)^{A_i} \prod_{j=1}^{\ell} \left( z - \frac{1}{2\varphi(d_j)}
  \right)^{B_j}, \qquad K \subset {\Bbb C} . $$
To this end, it suffices to show that
$$  p_n^*(z) \uto c\ p_0(z) = c (z^{A+B} + \lambda_1 z^{A+B-1} + \cdots +
    \lambda_{A+B}), \tag 3.6 $$
where
$$   c = \lim_{n \to \infty} \lambda_n^* = \left( \sum_{k=0}^{A+B}
     |\lambda_k| \right)^{-1}, \qquad \lambda_0=1. \tag 3.7 $$
Since $\{p_n^* \}$, $n \in  {\Bbb Z}_+$, is contained in ${\Bbb P}_{A+B}$ and
the sum of the
moduli of the coefficients of $p_n^*$, for each $n \in  {\Bbb Z}_+$, is
equal to one,
this family of polynomials is normal. Therefore \thetag{3.6} is obtained if we
prove that for all $\Lambda  \subset  {\Bbb Z}_+$ such that
$$  \lim_{n \to \infty,\ n \in \Lambda} p_n^*(z) = p_\Lambda  \tag 3.8 $$
then $p_\Lambda(z) = c\ p_0(z)$, where $p_0 (z)$ and $c$ are defined as above.
On the other hand, since $p_\Lambda  \in  {\Bbb P}_{A+B}$ and $p_\Lambda
\not\equiv 0$,
it is uniquely determined if we find its zeros and leading coefficient.
Finally, since the leading coefficient of $p_\Lambda$ is positive and the sum
of the moduli of its coefficients is one, its leading coefficient is
uniquely determined by its zeros. Therefore $p_\Lambda (z) = c\ p_0(z)$ if and
only if it is divisible by $p_0(z)$.

Because of the factor $S$ in $\Omega_n$ and since all the zeros of $L_{n-B}$
concentrate on $S_\mu$ we immediately obtain the following $A$ equations
$$  \multline
0 = \sum_{k=0}^{A+B} \lambda_n^* \lambda_{n,k}
   \left( \frac{L_{n+A-k}}{L_{n-B}} \right)^{(\nu)}(c_i),  \\
    i=1,\ldots,h;\ \nu=0,\ldots,A_i-1;\ n \geq n_0.
    \endmultline    \tag 3.9 $$
From \thetag{2.3} it follows that
$$  \left( \frac{L_{n+A-k}(z)}{L_{n-B}(z)} \right)^{(\nu)} \uto
   \left( \left( \frac{\varphi(z)}{2} \right)^{A+B-k} \right)^{(\nu)},
   \qquad K \subset {\Bbb C} \setminus S_\mu . \tag 3.10 $$
Relations \thetag{3.8}, \thetag{3.9} and \thetag{3.10}, together with the fact
that $\varphi$ holomorphic with $\varphi' \ne 0$ in ${\Bbb C} \setminus [-1,1]$,
imply (use induction on $\nu$)
$$  p_\Lambda^{(\nu)}\left( \frac{\varphi(c_i)}{2} \right) = 0,
  \qquad i=1,\ldots,h;\ \nu=0,\ldots,A_i-1 . \tag 3.11 $$
Take  $p(x) = T(x) L_{n-B}(x)/(x-d_j)^{\nu}$, $j=1,\ldots,\ell;\
\nu=1,\ldots,B_j$, in \thetag{3.1}.
Using \thetag{3.5} and multiplying by  $(\nu-1)! \lambda_n^* \tau_{n-B}^2$
we have the additional relations
$$  \multline
0 = \sum_{k=0}^{A+B} \lambda_n^* \lambda_{n,k}
   \frac{(\nu-1)! \tau_{n-B}}{\tau_{n+A-k}} \int \frac{l_{n+A-k}(x)
      l_{n-B}(x)}{(x-d_j)^\nu} \, d\mu(x), \\
       j=1,\ldots,\ell;\
      \nu=1,\ldots,B_j.
     \endmultline \tag 3.12 $$
From \thetag{2.4} and \thetag{2.12} it follows that
$$ \multline
(\nu-1)! \frac{\tau_{n-B}}{\tau_{n+A-k}}
   \int \frac{l_{n+A-k}(x)l_{n-B}(x)}{(x-z)^\nu} \, d\mu(x)
 \uto \left( \frac{-1}{(2\varphi(z))^{A+B-k} \sqrt{z^2-1}} \right)^{(\nu-1)},\\
    K \subset \overline{\Bbb C} \setminus S_\mu.
   \endmultline  \tag 3.13 $$
Relations \thetag{3.8}, \thetag{3.12} and \thetag{3.13}, together with the fact
that
$1/\varphi$ is holomorphic with  $(1/\varphi)'\ne 0 $ and $1/\sqrt{z^2-1} \ne
0$ in ${\Bbb C} \setminus [-1,1]$, give (use induction)
$$  p_\Lambda^{(\nu)}\left( \frac{1}{2\varphi(d_j)} \right) = 0,
    \qquad j=1,\ldots,\ell;\ \nu=0,\ldots,B_j-1. \tag 3.14 $$
From \thetag{3.11} and \thetag{3.14} follows that $p_\Lambda(z)$ is divisible
by $p_0(z)$. Therefore \thetag{3.6}--\thetag{3.7} hold and
$$  p_n(z) \uto p_0(z), \qquad K \subset {\Bbb C}. \tag 3.15 $$
From the definitions of $p_n$, $\Omega_n$, \thetag{3.10} with $\nu=0$ and
\thetag{3.15} we obtain
$$  \frac{R_n(z)}{L_{n-B}(z)} \uto p_0 \left( \frac{\varphi(z)}2 \right),
   \qquad K \subset {\Bbb C} \setminus S_\mu, $$
from which \thetag{3.2}, for $\nu  = 0$, immediately follows.

We continue by induction. Assume that \thetag{3.2} is true for $\nu  = k$. Note
that
$$  \frac{L_n^{(k)}}{L_n^{(k+1)}} \left( \frac{Q_n^{(k)}}{L_n^{(k)}} \right)'
   + \frac{Q_n^{(k)}}{L_n^{(k)}} = \frac{Q_n^{(k+1)}}{L_n^{(k+1)}} . $$
Hence, from \thetag{2.9} and \thetag{3.2} for $\nu =k$ we obtain
\thetag{3.2} for $\nu =k+1$.

Using Hurwitz' theorem (see \cite{1, p.~178}) \thetag{3.2}
implies that the zeros of $Q_n^{(\nu)}$,  for each fixed
$\nu \in {\Bbb Z}_+$, can only accumulate on $S_\mu$ as $n \to  \infty$.
Therefore, for sufficiently large $n$, the left hand members of
\thetag{3.3} and \thetag{3.4}
are meaningful (finite) on $K \subset \overline{\Bbb C} \setminus S_\mu$. Now
\thetag{3.3}
and \thetag{3.4} are trivial consequences of \thetag{3.2} combined respectively
with \thetag{2.5} and \thetag{2.9}. \qed
\enddemo

{\bf Remark 1.}
If $r(x)$ has constant sign on ${\Bbb R}$ then \thetag{3.3} and \thetag{3.4}
reduce to
\thetag{2.5} and \thetag{2.9} since then $r(x)\,d\mu(x)$ is a measure
of constant sign on its support. So only the case when $r$ has zeros
or poles in a non-symmetric way with respect to ${\Bbb R}$ is of interest.
\medskip

From \thetag{2.4}, \thetag{2.8} and \thetag{3.2} we can conclude that
$$  | Q_n^{(\nu)}(z) |^{1/n} \uto \frac{|\varphi(z)|}{2}, \qquad
    K \subset {\Bbb C} \setminus S_\mu. $$

\proclaim{Theorem 2}
Let $\mu  \in  \M(0,1)$. Then for all sufficiently large $n$
$(n\geq n_0)$
$$  \frac{1}{\kappa_n^2} := \int Q_n^2(x) r(x) \, d\mu(x) > 0 . \tag 3.16 $$
Moreover, for such $n$ the polynomials $Q_n$ satisfy a three-term recurrence
relation
$$    Q_{n+1}(z) = (z-\beta_n) Q_n(z) - \alpha_n^2 Q_{n-1}(z), \tag 3.17 $$
where
$$  \beta_n = \kappa_n^2 \int x Q_n^2(x) r(x) \, d\mu(x), \qquad
    \lim_{n \to \infty} \beta_n = 0 , \tag 3.18 $$
and
$$   \alpha_n^2 = \frac{\kappa_{n-1}^2}{\kappa_n^2}, \qquad
    \lim_{n \to \infty} \alpha_n^2 = \frac14 . \tag 3.19 $$
Finally
$$  \lim_{n \to \infty} \frac{\kappa_n^2}{\tau_n^2} =
    (-2)^{A-B} \frac{\prod_{j=1}^\ell \varphi^{B_j}(d_j) }
        {\prod_{i=1}^h \varphi^{A_i}(c_i) } . \tag 3.20 $$
\endproclaim

\demo{Proof}
If \thetag{3.16} were not true, then for an infinite set of indices $\Lambda $
one would have
$Q_n = Q_{n+1}$ for $n \in  \Lambda $, which contradicts the fact that for all
large $n$ the degree of $Q_n$ is $n$.

To prove \thetag{3.17} we observe that the right-hand member of that
relation is orthogonal with respect to $r(x)\,d\mu x)$ to all polynomials of
degree $\leq n-2$. If $n$ is large enough, $\beta_n$ and $\alpha_n^2$ as
given by the left hand of \thetag{3.18} and \thetag{3.19} are well defined
(finite) and furthermore the degree of $Q_{n-1}$ is $n-1$ the degree of
$Q_n$ is $n$. We also have
$$ \int \left[ (x-\beta_n) Q_n(x) - \alpha_n^2 Q_{n-1}(x) \right]
   Q_{n-1}(x) r(x) \, d\mu(x) = \frac{1}{\kappa_n^2} -
   \frac{\alpha_n^2}{\kappa_{n-1}^2} = 0, $$
and
$$ \multline
\int \left[ (x-\beta_n) Q_n(x) - \alpha_n^2 Q_{n-1}(x) \right]
   Q_n(x) r(x) \, d\mu(x)  \\
   = \int x Q_n^2(x) r(x) \, d\mu(x) - \frac{\beta_n}{\kappa_n^2} = 0.
 \endmultline   $$
Since the degree of $(z-\beta_n)Q_n(z)-\alpha_n^2 Q_{n-1}(z)$ is $n+1$
(the degree of $Q_{n+1}$ for large $n$), we see that \thetag{3.17} follows
from the above.

To prove the limit relations in \thetag{3.18} and \thetag{3.19} we proceed as
follows. Equation \thetag{3.17} may be written as
$$   z Q_n(z) = Q_{n+1}(z) + \beta_n Q_n(z) + \alpha_n^2 Q_{n-1}(z) . $$
Multiplying either side by $S(z)$ and using \thetag{3.5} we have
$$ \multline
\sum_{k=0}^{A+B} \lambda_{n,k} z L_{n+A-k}(z)
   = \sum_{k=0}^{A+B} \lambda_{n+1,k} L_{n+A+1-k}(z)  \\
   + \sum_{k=0}^{A+B} \beta_n \lambda_{n,k} L_{n+A-k}(z)
   + \sum_{k=0}^{A+B} \alpha_n^2 \lambda_{n-1,k} L_{n+A-1-k}(z).
   \endmultline \tag 3.21 $$
Using the three-term recurrence relation for the polynomials $L_n$
$$   z L_n(z) = L_{n+1}(z) + b_n L_n(z) + a_n^2 L_{n-1}(z), $$
the left hand of \thetag{3.21} can be written as
$$ \multline
\sum_{k=0}^{A+B} \lambda_{n,k} z L_{n+A-k}(z)
   = \sum_{k=0}^{A+B} \lambda_{n,k} L_{n+A+1-k}(z)  \\
   + \sum_{k=0}^{A+B} b_{n+A-k} \lambda_{n,k} L_{n+A-k}(z)
   + \sum_{k=0}^{A+B} a_{n+A-k}^2 \lambda_{n,k} L_{n+A-1-k}(z).
   \endmultline \tag 3.22 $$
Since $\{L_n\}$, $n\in {\Bbb Z}_+$, is a basis, comparing the coefficients
corresponding to $L_{n+A}$ and $L_{n+A-1}$ in the right-hand members of
\thetag{3.21} and \thetag{3.22} gives
$$ \aligned
\lambda_{n+1,1} + \beta_n \lambda_{n,0} &= \lambda_{n,1} + b_{n+A}
\lambda_{n,0},  \\
\lambda_{n+1,2} + \beta_n \lambda_{n,1} + \alpha_n^2 \lambda_{n-1,0} &=
\lambda_{n,2} + b_{n+A-1} \lambda_{n,1} + a_{n+A}^2 \lambda_{n,0}.
  \endaligned  \tag 3.23 $$
We know from \thetag{3.15} that
$$  \lambda_{n,0} = \lambda_{n-1,0} = 1, \quad
    \lim_{n \to \infty} \lambda_{n+1,1} = \lim_{n \to \infty} \lambda_{n,1},
    \quad
    \lim_{n \to \infty} \lambda_{n+1,2} = \lim_{n \to \infty} \lambda_{n,2}, $$
and  since $\mu \in \M(0,1)$ we also have
$$  \lim_{n \to \infty} b_n = 0, \quad \lim_{n \to \infty} a_n^2 = \frac14. $$
From all this and equations \thetag{3.23} we get \thetag{3.18} and
\thetag{3.19}.

In  order  to  prove  \thetag{3.20}, note  that  for  $n \geq  n_0$ the
degree of $Q_n-TL_{n-B}$ is $\leq n-1$. Therefore, from this and
\thetag{3.5}
$$ \multline
  \frac{\tau_n^2}{\kappa_n^2} = \tau_n^2 \int Q_n^2(x) r(x) \, d\mu(x)
  = \tau_n^2 \int L_{n-B}(x) Q_n(x) S(x) \, d\mu(x) \\
  = \lambda_{n,A+B} \tau_n^2 \int L_{n-B}^2(x) \, d\mu(x)
  = \lambda_{n,A+B} \frac{\tau_n^2}{\tau_{n-B}^2} .
  \endmultline   \tag 3.24  $$
Since
$$  \lim_{n \to \infty} \lambda_{n,A+B} = \lambda_{A+B}
   = \left( \frac{-1}2 \right)^{A+B} \frac{ \prod_{i=1}^h
      \varphi^{A_i}(c_i) }{ \prod_{j=1}^\ell \varphi^{B_j}(d_j) }
      \tag 3.25 $$
(see \thetag{3.15}, \thetag{3.6} and the relation preceding \thetag{3.6}),
we can now use \thetag{2.4}, \thetag{3.24} and \thetag{3.25} to obtain
\thetag{3.20}. This completes the proof of the theorem. \qed
\enddemo

Writing
$$  \alpha_n = \left( \frac{ \int Q_n^2(x) r(x) \, d\mu(x) }
    { \int Q_{n-1}^2(x) r(x)\, d\mu(x) } \right)^{1/2}, $$
(cfr. \thetag{3.18} and \thetag{3.19}), we have by \thetag{3.19}
$$  \lim_{n \to \infty} \alpha_n = \frac12 . $$
Choosing $n_0$ such that
$$  \kappa_{n_0} = \left( \int Q_{n_0}^2(x) r(x) \, d\mu(x)
\right)^{-1/2} \ne 0 ,\  $$
we can define $\kappa_n$ for $n > n_0$ recursively by
$$  \kappa_{n+1} = \frac{\kappa_n}{\alpha_{n+1}} . $$
With this selection we take $q_n = \kappa_n Q_n$ and the
sequence $\{q_n\}$, $n \geq  n_0$ becomes orthonormal with respect to
$d\rho(x)=r(x)\,d\mu(x)$, with the degree of $q_n$ equal to $n$. Moreover,
because of \thetag{3.17} these polynomials satisfy the
three-term recurrence relation
$$  x q_n(x) = \alpha_{n+1} q_{n+1}(x) + \beta_n q_n(x) + \alpha_n q_{n-1}(x),
   \qquad n \geq n_0, \tag 3.26 $$
with
$$  \lim_{n \to \infty} \alpha_n = \frac12, \quad
  \lim_{n \to \infty} \beta_n = 0.    \tag 3.27 $$
Moreover, for all fixed $\nu  \in {\Bbb Z}_+$, the zeros of $q_n^{(\nu)}$,
$n \geq n_0$, can
only accumulate on $S_\mu$, as $n \to  \infty$, and
$$  \frac{q_{n+1}^{(\nu)}(z)}{q_n^{(\nu)}(z)} \uto \varphi(z),
   \qquad K \subset {\Bbb C} \setminus S_\mu, \tag 3.28 $$
together with
$$  \frac1n \frac{q_n^{(\nu+1)}(z)}{q_n^{(\nu)}(z)} \uto
\frac{1}{\sqrt{z^2-1}}, \qquad K \subset \overline{\Bbb C} \setminus S_\mu.
\tag 3.29 $$

From previous results it is easy to prove also that for all fixed
$\nu \in {\Bbb Z}$
$$   \int |q_n(x)q_{n+\nu}(x)|\ |d\rho(x)| \leq \C < \infty . \tag 3.30 $$
Here, as usual, $|d\rho|$ is the total variation of the complex measure
$\rho$ (see, e.g., \cite{23}).
In fact (see \thetag{3.5}) one has
\hfuzz=7pt
$$ \multline
  \int |q_n(x)q_{n+\nu}(x)|\ |d\rho(x)| =
  |\kappa_n \kappa_{n+\nu}| \int |Q_n(x)S(x) Q_{n+\nu}(x) S(x)|
   \frac{d\mu(x)}{|S(x)T(x)|} \\
   \leq |\kappa_n \kappa_{n+\nu}|
   \sum_{k=0}^{A+B} \sum_{m=0}^{A+B} \int |\lambda_{n,k} \lambda_{n+\nu,m}
   L_{n+A-k}(x)L_{n+\nu+A-m}(x)| \frac{d\mu(x)}{|S(x)T(x)|} \\
   = \frac{|\kappa_n \kappa_{n+\nu}|}{\tau_n \tau_{n+\nu}}
   \sum_{k=0}^{A+B} \sum_{m=0}^{A+B}
   \frac{\tau_n\tau_{n+\nu}|\lambda_{n,k}\lambda_{n+\nu,m}|}
   {\tau_{n+A-k}\tau_{n+A+\nu-m}}
   \int |l_{n+A-k}(x) l_{n+A+\nu-m}(x)| \frac{d\mu(x)}{|S(x)T(x)|} .
   \endmultline $$
\hfuzz=0.1pt
Set $\C = \inf  \{|S(x)T(x)|: x \in  S_\mu \}$ $(>0)$. Then using the
Cauchy-Schwarz inequality one has
$$ \int |q_n(x)q_{n+\nu}(x)|\ |d\rho(x)| \leq
  \frac{1}{\C}  \frac{|\kappa_n \kappa_{n+\nu}|}{\tau_n \tau_{n+\nu}}
  \sum_{k=0}^{A+B} \sum_{m=0}^{A+B}
  \frac{\tau_n}{\tau_{n+A-k}} \frac{\tau_{n+\nu}}{\tau_{n+A+\nu-m}}
  |\lambda_{n,k} \lambda_{n+\nu,m}| , $$
and \thetag{3.30} immediately follows using \thetag{2.4}, \thetag{3.15} and
\thetag{3.20}.

In \cite{14} A. Magnus introduced a more general class of complex measures
to
which, for large $n$, a sequence of (non-Hermitian) orthonormal polynomials may
be associated satisfying \thetag{3.26}---\thetag{3.30}. These are all complex
measures on
$[-1,1]$ of the form $d\rho(x) = g(x)\omega(x)\, dx$, where $g$ is a
non-vanishing
complex valued continuous function on $[-1,1]$ and $\omega$ is a positive
almost everywhere real integrable function on $[-1,1]$.

When $\beta_n \in {\Bbb R}$, $\alpha_n > 0$, the existence of a positive
measure $\rho$ supported in ${\Bbb R}$ with respect to which $\{q_n\}$ is
orthonormal
is guaranteed from the recursion formula by Favard's theorem. If, in
addition, \thetag{3.27} holds  $(\rho \in \M(0,1))$, then $\rho$ has the
structure
described in Section 1 and automatically \thetag{3.28}--\thetag{3.30}
follows
(\thetag{3.30} is trivial for positive measures), see \cite{14}, \cite{19},
\cite{24} and \cite{27}.

Sequences of polynomials satisfying general recurrence relations of type
\thetag{3.26}
have been studied. See, for example, \cite{4} for representation theorems and
\cite{7} for bounds on the zeros. Such relations are satisfied by the
denominators
of diagonal Pad\'e approximants and certain types of continued fractions,
hence their importance in approximation theory.
Such sequences of polynomials need not come from complex measures (for example,
they also arise from orthogonality relations of type \thetag{1.6}). In case
they do,
what restrictions on the class of measures, where we seek the solution of the
corresponding complex moment problem, determine a unique measure?
This problem is  commonly known as the question of determinacy.

Recently, one of the authors in a joint paper with E. Torrano and D. Barrios
\cite{13} studied the location of zeros of polynomials given by \thetag{3.26}
under the assumptions:
$$  \sup_n (|\alpha_n|,|\beta_n|) \leq C < \infty, \quad
   \lim_{n \to \infty} {\frak Im\ } \alpha_n = \lim_{n \to \infty} {\frak Im\ }
 \beta_n = 0 .  $$
In particular, from their results it follows that if \thetag{3.27} holds,
then automatically \thetag{3.28} and \thetag{3.29} hold uniformly on compact
subsets
of ${\Bbb C} \setminus ([-1,1] \cup E)$, where $E$ is at most a denumerable
set of isolated points in  ${\Bbb C} \setminus [-1,1]$
for which $E'\subset [-1,1]$. Moreover, in this case
\thetag{3.27} also implies uniqueness of the generating measure (should it
exist)
with support contained in $[-1,1]\cup E$. An open problem is still whether under
these circumstances \thetag{3.30} is satisfied.

\proclaim{Definition 1}
We say that $\rho \in \MC(0,1)$ if the corresponding
orthonormal polynomials satisfy \thetag{3.26}, \thetag{3.27} and \thetag{3.30},
where $S_\rho =[-1,1]\cup E$, and  $E$  is  at  most a denumerable set of
isolated points in ${\Bbb C} \setminus [-1,1]$.
\endproclaim

\proclaim{Theorem 3}
Let $\rho  \in  \MC(0,1)$. Then for all $f$ continuous on
$S_\rho$ and fixed $\nu \in {\Bbb Z}$
$$  \lim_{n \to \infty} \int f(x) q_n(x) q_{n+\nu}(x) \, d\rho(x)
    = \frac{1}{\pi} \int_{-1}^1 \frac{f(x) t_{|\nu|}(x) \, dx}{\sqrt{1-x^2}},
  \tag 3.31 $$
where $t_{|\nu|}$ is the $|\nu|$th Chebyshev orthonormal
polynomial of the first kind.
\endproclaim

\demo{Proof}
If $f(x)$ is a power of $x$, then the proof is carried out in the same fashion
as Theorem 4.2.13 in \cite{19}. Therefore, \thetag{3.31} holds
for all polynomials. If $f(x)$ is a continuous function on $S_\rho$ then,
since  $S_\rho^c$ is connected and $S_\rho^o = \emptyset$,
according to a theorem of M. Lavrentiev \cite{6, Thm.~8.7 on p.~48}
 $f$ may be uniformly approximated on
$S_\rho$ by polynomials. This combined with \thetag{3.30} implies
\thetag{3.31} for $f$. \qed
\enddemo

\proclaim{Corollary 1}
Assume that $\rho  \in  \MC(0,1)$. Then for all
$\nu, \eta  \in  {\Bbb Z}_+$ with $\nu \geq \eta$
$$      \frac{q_n^{(\nu)}(z)}{n^{\nu-\eta} q_n^{(\eta)}(z)}
    \uto \left( \frac{1}{\sqrt{z^2-1}} \right)^{\nu-\eta},
    \qquad K \subset \overline{\Bbb C} \setminus S_\rho, \tag 3.32 $$
and if  $m,k \in {\Bbb Z}$ with $k \leq m$ and $\nu, \eta \in {\Bbb Z}_+$
$$ \multline
   \frac{(-1)^{\eta-1}(\eta-1)!}{n^\nu}
     \int q_{n-m}^{(\nu)}(x)q_{n-k}(x) \frac{d\rho(x)}{(z-x)^\eta} \\
     \uto
     \left[ \left( \frac{1}{\sqrt{z^2-1}} \right)^\nu
     \frac1\pi \int \frac{t_{m-k}(x) \, dx}{(z-x) \sqrt{1-x^2}}
     \right]^{(\eta-1)}, \qquad K \subset \overline{\Bbb C} \setminus S_\rho.
     \endmultline \tag 3.33 $$
\endproclaim

\demo{Proof}
The asymptotic behavior \thetag{3.32} follows from \thetag{3.29} and the fact
that
$$  \frac{q_n^{(\nu)}(z)}{n^{\nu-\eta} q_n^{(\eta)}(z)}
   = \prod_{j=\eta}^{\nu-1}  \frac{q_n^{(j+1)}(z)}{n q_n^{(j)}(z)} . $$
Let us prove \thetag{3.33} for $\eta=1$. Since $k \leq m$ it follows that
$$  \int \frac{q_{n-m}^{(\nu)}(z)-q_{n-m}^{(\nu)}(x)}{z-x} \, q_{n-k}(x)
   \, d\rho(x) = 0 . $$
Therefore
$$  \aligned
\frac{1}{n^\nu}
 \int \frac{q_{n-m}^{(\nu)}(x)q_{n-k}(x)}{z-x} \, d\rho(x) &=
 \frac{q_{n-m}^{(\nu)}(z)}{n^\nu} \int \frac{q_{n-k}(x)}{z-x} \, d\rho(x) \\
 &= \frac{q_{n-m}^{(\nu)}(z)}{n^\nu q_{n-m}(z)}
  \int \frac{q_{n-m}(x)q_{n-k}(x)}{z-x} \, d\rho(x) .
  \endaligned \tag 3.34 $$
From \thetag{3.30} and \thetag{3.32} we have that, given
$K \subset \overline{\Bbb C} \setminus S_\rho$,
this family of
functions is uniformly bounded for $z \in  K$. Thus
\thetag{3.33} with $\eta  = 1$
 follows from pointwise convergence. This is guaranteed by \thetag{3.31} with
 $f(x) = 1/(z-x)$, \thetag{3.32} and \thetag{3.34}.

For arbitrary $\eta  \geq  1$ we observe that
$$ \frac{1}{n^\nu}
\int \frac{q_{n-m}^{(\nu)}(x)q_{n-k}(x)}{(z-x)^\eta} \, d\rho(x) =
\frac{(-1)^{\eta-1}}{(\eta-1)!} \left(
 \frac{1}{n^\nu} \int \frac{q_{n-m}^{(\nu)}(x)q_{n-k}(x)}{z-x} \, d\rho(x)
 \right)^{(\eta-1)}, $$
and we obtain \thetag{3.33} using that formula with $\eta =1$.
Formula \thetag{2.13} allows to  express \thetag{3.33} in terms
of $\varphi$. \qed
\enddemo

We conclude this section with two lemmas.

\proclaim{Lemma 3}
 Suppose that $\rho  \in  \MC(0,1)$, $\eta ,\nu \in {\Bbb Z}_+$. Then
$$  \left( \frac{1}{q_n^{(\eta)}(z)} \right)^{(\nu)}
   = \sum_{k=0}^\nu g_{n,k}  \frac{n^k}{q_n^{(\eta)}(z)}, \tag 3.35 $$
where $g_{n,k}$, $k = 0,\ldots,\nu$, are rational functions whose poles
accumulate on $S_\rho$ as $n \to  \infty$. Moreover
$$  g_{n,k} \uto g_k \in H(\overline{\Bbb C} \setminus S_\rho),
  \qquad K \subset \overline{C} \setminus S_\rho . \tag 3.36 $$
In particular, for each $\eta  \in  {\Bbb Z}_+$
$$   g_{n,\nu} \uto \left( \frac{-1}{\sqrt{z^2-1}} \right)^\nu,
   \qquad K \subset \overline{C} \setminus S_\rho. $$
\endproclaim

\demo{Proof}
 We proceed by induction. Fix $\eta  \in  {\Bbb Z}_+$; for $\nu  = 0$
the equations \thetag{3.35}, \thetag{3.36} and the statement
about $g_{n,0}$ follow easily by taking $g_{n,0} = 1$.
Assume the lemma holds for $\nu =m$, and let us prove it for $\nu =m+1$.
Using the induction hypothesis we obtain
$$ \align
 \left( \frac{1}{q_n^{(\eta)}(z)} \right)^{(m+1)} &=
 \left( \sum_{k=0}^m g_{n,k}(z) \frac{n^k}{q_n^{(\eta)}(z)}
  \right)' \\
  &= \sum_{k=0}^m g_{n,k}'(z) \frac{n^k}{q_n^{(\eta)}(z)}
    - \sum_{k=0}^m g_{n,k} \frac{q_n^{(\eta+1)}(z)}{n q_n^{(\eta)}(z)}
      \frac{n^{k+1}}{q_n^{(\eta)}(z)} \\
  &= g_{n,0}'(z) \frac{1}{q_n^{(\eta)}(z)}
    + \sum_{k=1}^m \left( g_{n,k}'(z) - g_{n,k-1}
    \frac{q_n^{(\eta+1)}(z)}{n q_n^{(\eta)}(z)} \right)
    \frac{n^k}{q_n(z)} \\
    &\quad -\ g_{n,m}(z) \frac{q_n^{(\eta+1)}(z)}{n q_n^{(\eta)}(z)}
    \frac{n^{m+1}}{q_n^{(\eta)}(z)} .
  \endalign $$
Using the induction hypothesis together with \thetag{3.32},
the  lemma readily follows for $\nu =m+1$.  \qed
\enddemo

A repetition of the arguments in the proofs of Theorems 1, 2 and the
consequences thereafter gives:

\proclaim{Lemma 4}
 Let $\rho \in \MC(0,1)$ and $r$ be a rational function whose
zeros and poles lie in ${\Bbb C} \setminus S_\rho$. Then $r\,d\rho \in
\MC(0,1)$.
 Moreover, the orthogonal polynomials with respect to
$r\,d\rho$ and $d\rho$ have relative asymptotic behavior as described by
\thetag{3.2} and \thetag{3.20}.
\endproclaim

\head 4. Relative asymptotics for generalized discrete Sobolev inner
    products \endhead

\proclaim{Definition 2}
Let $\mu$ be a complex measure with compact support $S_\mu \subset {\Bbb C}$
and $c_1,c_2,\ldots,c_m \in {\Bbb C}$. A generalized discrete Sobolev
inner product is an expression of the form
$$  \langle h,g \rangle = \int hg\, d\mu
     + \sum_{j=1}^m \sum_{i=0}^{N_j} h^{(i)}(c_j) \L_{j,i}(g;c_j),
     \tag 4.1 $$
where $\L_{j,i}(g;c_j)$ is the evaluation at $c_j$ of the linear ordinary
differential operator $\L_{j,i}$ with constant coefficients acting on
$g$ and $\L_{j,N_j} \not\equiv 0$, $j=1,\ldots,m$.
\endproclaim

Since the operators in \thetag{4.1} are evaluated at a single point $c_j$, the
assumption that they have constant coefficients is not a restriction, but the
results hereafter are better described this way. The linearity is essential.

Denote by $S_n$, $n \in {\Bbb Z}_+$, the monic polynomial of the least
degree, such that
$$  \langle p,S_n \rangle = 0, \qquad p \in {\Bbb P}_{n-1}. \tag 4.2 $$
The existence of $S_n$, for each $n \in {\Bbb Z}_+$, follows from solving
a linear system of $n$ homogeneous equations of $n+1$ unknowns. Uniqueness
readily follows from being of the least degree. If the inner product is positive
definite, then the degree of $S_n$ is $n$ and thus all the $S_n$'s are distinct.
In general this is not so and different $n$'s may have the same $S_n$.

An inner product of the form \thetag{4.1} may not be symmetric. Therefore
we must distinguish between right-orthogonal polynomials defined by
\thetag{4.2} and left-orthogonal polynomials, which are defined by
$\langle S_n,p \rangle = 0,\ p \in {\Bbb P}_{n-1}$. However, since
$$  \langle h,g \rangle =
   \int gh\, d\mu + \sum_{j=1}^m \sum_{i=0}^{N_j^*} g^{(i)}(c_j)
      \L_{j,i}^*(h;c_j) = \langle g,h \rangle^*, $$
it follows that left-orthogonal polynomials with respect to $\langle
\cdot,\cdot \rangle$ are right-or\-tho\-go\-nal polynomials with respect to
$\langle \cdot,\cdot \rangle^*$. Therefore it suffices to study right-orthogonal
polynomials, which we will simply call orthogonal polynomials.

Special cases of generalized discrete Sobolev inner products are
\thetag{1.1} and the inner product associated with \thetag{1.6}, that is
$$  \int hg\, d\mu + \sum_{j=1}^m \sum_{i=0}^{N_j}
    h^{(i)}(c_j) \sum_{k=i}^{N_j} A_{j,k} {k \choose i}
     g^{(k-i)}(c_j) . \tag 4.3 $$

Given $j=1,\ldots,m$, let $J_j$ be the maximum order of the differential
operator $\L_{j,i}$, $i=0,\ldots,N_j$. Then
$$  \L_{j,i}(g) = \sum_{k=0}^{J_j} \gamma_{i,k}^j g^{(k)} . $$
Let $\Gamma_j = (\gamma_{i,k}^j)$, $i=0,\ldots,N_j$, $k=0,\ldots,J_j$, be
the matrix of the coefficients of the $\L_{j,i}$. Denote by
$\Gamma_j^*$ the matrix obtained from $\Gamma_j$ after deleting all the
rows and columns with zero vectors.

\proclaim{Definition 3}
We say that $\langle \cdot,\cdot \rangle$ is a regular generalized
discrete Sobolev inner product if for each $j=1,\ldots,m$ the matrix
$\Gamma_j^*$ is a square matrix with determinant different from zero.
We denote by $I_j$ the dimension of $\Gamma_j^*$.
\endproclaim

This definition means that the total number of different derivatives appearing
in the $\L_{j,i}$'s equals for each $j$ the number of these operators not
identically equal to zero, plus the fact that the non-zero operators be
linearly independent. Obviously this is the case with \thetag{1.1}
where for each $j$ the matrix $\Gamma_j^*$ is a diagonal matrix, and with
\thetag{4.3} where the corresponding $\Gamma_j^*$ are triangular square
matrices.
Both of these inner products are symmetric (for \thetag{4.3} this follows
easily from its equivalent form \thetag{1.6}), but non-symmetric
regular generalized discrete Sobolev inner products are easy to construct.
Obviously $\langle \cdot,\cdot \rangle$ is symmetric if for each
$j$ the matrix $\Gamma_j^*$ is symmetric. From this fact it immediately
follows that not all such symmetric inner products are regular.

The theorem which we are about to prove can easily be extended to certain
non-regular inner products of type \thetag{4.1}, but we restrict our attention
to regular inner products for three reasons. Firstly, they contain the most
interesting cases \thetag{1.1} in its complex version $(M_{j,i} \in {\Bbb C},
M_{j,N_j} \neq 0, c_j \in {\Bbb C})$ and \thetag{4.3}. Secondly,
the notation and technicalities rapidly increase with generality, making the
reading (and writing) difficult. Thirdly, because we are sure that the
best statement we can prove does not have a final character. The proof
of Theorem 4 contains the main ingredients of our approach in solving the
general problem.

\proclaim{Theorem 4}
Consider a regular inner product of type \thetag{4.1} such that $\mu \in
\MC(0,1)$ and $c_1,\ldots,c_m \in {\Bbb C} \setminus S_\mu$.
Let $\{ L_n \}$, $n \in {\Bbb Z}_+$, be the sequence of monic orthogonal
polynomials with respect to $\mu$ and $\{ S_n \}$, $n \in {\Bbb Z}_+$, the
monic orthogonal
polynomials with respect to the given inner product. Then, for all
sufficiently large $n$, the degree of $S_n$ is $n$ and each point $c_j$
attracts exactly $I_j$ zeros of $S_n$, while the rest of the zeros
concentrate on $S_\mu$. Also, for each fixed $\nu \in {\Bbb Z}_+$
$$  \frac{S_n^{(\nu)}(z)}{L_n^{(\nu)}(z)} \uto
  \prod_{j=1}^m \left( \frac{(\varphi(z)-\varphi(c_j))^2}{2\varphi(z)(z-c_j)}
   \right)^{I_j}, \qquad K \subset \overline{\Bbb C} \setminus S_\mu .
   \tag 4.4 $$
\endproclaim

\demo{Proof}
Write
$$  s(z) = \prod_{j=1}^m (z-c_j)^{N_j+1}, \quad
    \tilde{s}_j(z) = \frac{s(z)}{(z-c_j)^{N_j+1}}.  $$
Taking $p(z) = s(z) p_1(z)$ in \thetag{4.2}, where $p_1$ is an arbitrary
polynomial of degree $\leq n-A-1$ with $A$ the degree of $s$, we have
$$  0 = \langle sp_1 , S_n \rangle = \int p_1 S_n \, d\rho, \tag 4.5 $$
where $d\rho = s\, d\mu$. By Lemma 4 we know that $\rho \in \MC(0,1)$,
and thus, by the same lemma, we also know the asymptotic behavior of the monic
orthogonal polynomials $Q_n$, $n \in {\Bbb Z_+}$, for the measure
$\rho$ relative to the
polynomials $L_n$, $n \in {\Bbb Z_+}$. Therefore, instead of
\thetag{4.4} we shall study the behavior of the ratio
$S_n^{(\nu)}/Q_n^{(\nu)}$, $n \in {\Bbb Z}_+$. Formula \thetag{4.4}
for arbitrary $\nu \in {\Bbb Z}_+$ follows, as before, by induction on $\nu$.
Thus we will restrict our attention to the case $\nu=0$. Since the degree
of $Q_n$ is $n$ for all sufficiently large $n$, which in turn implies that
for such $n$'s
$$  \int Q_n^2 \, d\rho \neq 0, $$
we can conclude that any polynomial $S_n$ satisfying  \thetag{4.5}
must be of degree at least $n-A$ for $n\geq n_0$. Now it is easy to see
that
$$   S_n(z) = \sum_{k=0}^A \lambda_{n,k} Q_{n-k}(z), \tag 4.6 $$
since any such polynomial satisfies \thetag{4.5} and the $A+1$ parameters
$\lambda_{n,k}$, $k=0,\ldots,A$ are sufficient to guarantee a nontrivial
solution of the remaining $A$ homogeneous linear relations
$$  0 = \langle t^\nu , \sum_{k=0}^A \lambda_{n,k} Q_{n-k} \rangle ,
    \qquad \nu=0,\ldots,A-1, $$
which determine $S_n$.

As in Theorem 1, the proof of the corresponding asymptotic formula for
the sequence
$\{ S_n/Q_n \}$, $n \in {\Bbb Z }_+$, relies on proving that for each
$k=0,\ldots,A,$ one has the asymptotic behavior $\lambda_{n,k} \to \lambda_k$
and finding the limits $\lambda_k$. These parameters are not the same
as those appearing in \thetag{3.5}, but in the proof they play a similar role.
Write
$$  \Omega_n(z) := \frac{S_n(z)}{Q_{n-A}(z)}
    = \sum_{k=0}^A \lambda_{n,k} \frac{Q_{n-k}(z)}{Q_{n-A}(z)}. $$
Set
$$  \lambda_n^* = \left( \sum_{k=0}^A |\lambda_{n,k}| \right)^{-1}
  \quad ( < \infty), $$
and
$$  h_n(z) = \sum_{k=0}^A \lambda_{n,k} z^{A-k}, \quad
    h_n^*(z) = \lambda_n^* h_n(z). $$
We shall prove that
$$  h_n(z) \uto h_0(z) :=
   \prod_{j=1}^m \left( z- \frac{\varphi(c_j)}{2} \right)^{I_j}
      \left( z - \frac{1}{2\varphi(c_j)} \right)^{N_j+1-I_j},
       \qquad K \subset {\Bbb C}. \tag 4.7 $$
As in Section 3, this follows by proving that
$$  h_n^*(z) \uto c\ h_0(z) = c ( z^A + \lambda_1 z^{A-1} + \cdots +
\lambda_A ), \tag 4.8 $$
where
$$  c = \lim_{n \to \infty}  \lambda_n^*
      = \left( \sum_{k=0}^A |\lambda_k| \right)^{-1} , \qquad
      \lambda_0 = 1. $$
Since the sum of the moduli of the coefficients of $h_n^*$ is one, the family
$\{ h_n^* \}$, $n \in {\Bbb Z}_+$, is normal and \thetag{4.8} holds if
we prove that for each subset $\Lambda \subset {\Bbb Z}_+$ for which
$$  \lim_{n \to \infty,\ n \in \Lambda} h_n^*(z) = h_\Lambda(z),
         \tag 4.9 $$
the limit $h_\Lambda$ is divisible by $h_0$. To this end, let us find a
convenient system of $A$ linear equations for the parameters
$\lambda_{n,k}$, $k=0,\ldots,A$.
Fix $j \in \{1,2,\ldots,m\}$ and put $p(z)=\tilde{s}_j(z) (z-c_j)^I$,
$I=0,\ldots,N_j$, in \thetag{4.2}. Then
$$   0 = \int S_n(x) \frac{d\rho(x)}{(x-c_j)^{N_j+1-I}}
    + \sum_{i=I}^{N_j} I! {i \choose I} \tilde{s}_j^{(i-I)}(c_j)
    \L_{j,i}(S_n;c_j). \tag 4.10 $$
Since $j$ is fixed we will drop it from the notation in the rest of the proof.
Thus, in what follows we use $N_j=N$, $\L_{j,i}=\L_i$, $c_j=c$, and
$\tilde{s}_j=s$ (which differs from the $s$ used higher). Note that
$s(c)=\tilde{s}_j(c_j) \neq 0$. This along with \thetag{4.10} allows us to prove
that there exist polynomials $p_i=p_{j,i}$ independent of $n$ such that
the degree of $p_i$ is $N+1-i$, $p_i(0)=0$ and
$$  0 = \int S_n(x) p_i\left(\frac1{x-c}\right) \, d\rho(x)
  + \L_i(S_n;c), \qquad i=0,\ldots, N. \tag 4.11 $$
Indeed, system \thetag{4.11} is equivalent to \thetag{4.10} and is obtained
from it in view of the triangular structure of the sum in \thetag{4.10}.

We will combine the equations in \thetag{4.11} in two different groups: those
corresponding to the operators $\L_i \neq 0$ and those where $\L_i \equiv 0$.
Let $i_1,\ldots,i_I$ be those indices $i$ for which $\L_i \not\equiv 0$.
According to the assumption of regularity, $I=I_j$ is equal to the dimension of
$\Gamma_j^*$ (see Definition 3). Let $d_1,\ldots,d_I$ be the orders of the
derivatives appearing in $\L_{i_1},\ldots,\L_{i_I}$ and recall that
$\Gamma_j^*$ is a square matrix. Since $\det(\Gamma_j^*) \neq 0$ we can
further transform the equations in \thetag{4.11} corresponding to the indices
$i_1,\ldots,i_I$ and find polynomials $p_{i_1}^*,\ldots,p_{i_I}^*$,
independent of $n$, such that the degree of $p_{i_r}^*$ is $\leq N+1$,
$p_{i_r}^*(0)=0$, $r=1,\ldots,I$, and
$$  0 = \int S_n(x) p_{i_r}^*\left( \frac{1}{x-c} \right) \, d\rho(x)
    + S_n^{(d_r)}(c), \qquad r=1,\ldots,I. \tag 4.12 $$
In this group of equations the main role is played by the second term on
the right hand side. In order to estimate the integral term, we proceed
as follows. Take $q_n$ the $n$th orthonormal polynomial with respect to $\rho$.
Since for each $\ell \in {\Bbb N}$ (the set of non-negative integers)
$$  \left( q_{n-A}(c) - \sum_{\eta=0}^{\ell-1}
   \frac{q_{n-A}^{(\eta)}(x)}{\eta !} (c-x)^\eta \right) (x-c)^{-\ell} $$
is a polynomial (in $x$ and $c$) of degree $\leq n-A-\ell < n-k$, it follows
by orthogonality that
$$  q_{n-A}(c) \int \frac{q_{n-k}(x) \, d\rho(x)}{(x-c)^\ell}
   = \sum_{\eta=0}^{\ell-1} (-1)^\eta
     \int \frac{q_{n-A}^{(\eta)}(x) q_{n-k}(x)}{\eta ! (x-c)^{\ell-\eta}}
        \, d\rho(x). \tag 4.13 $$
Since the degree of $p_{i_r}^*$ is $\leq N+1$, we obtain from
\thetag{4.13} and \thetag{3.33}
$$   \int q_{n-k}(x) p_{i_r}^* \left( \frac{1}{x-c} \right) \, d\rho(x)
    = O \left( \frac{n^N}{q_{n-A}(c)} \right) . \tag 4.14 $$
In what follows, each $O(\cdot)$ may be a different function
of $n$.

Now let us turn to the second term in \thetag{4.12}. We may assume that the
indices $i_r$ are taken so that $0 \leq d_1 < \cdots < d_I$. Multiply
both sides of \thetag{4.12} by $\lambda_n^* \kappa_{n-A}/q_{n-A}^{(d_r)}(c)$,
where $\kappa_n$ is the leading coefficient of the $n$th orthogonal polynomial
$q_n$ for the measure $\rho$. Using \thetag{4.6} and \thetag{4.14} we have
$$  \aligned
   0 &= \sum_{k=0}^A \lambda_n^* \lambda_{n,k} \frac{\kappa_{n-A}}{\kappa_{n-k}}
     \left[  \int \frac{q_{n-k}(x)}{q_{n-A}^{(d_r)}(c)}
     p_{i_r}^*\left( \frac{1}{x-c} \right) \, d\rho(x)
     + \frac{q_{n-k}^{(d_r)}(c)}{q_{n-A}^{(d_r)}(c)} \right] \\
  &= \sum_{k=0}^A \lambda_n^* \lambda_{n,k} \frac{\kappa_{n-A}}{\kappa_{n-k}}
     \left[ O \left( \frac{n^N}{q_{n-A}(c) q_{n-A}^{(d_r)}(c)} \right)
     + \frac{q_{n-k}^{(d_r)}(c)}{q_{n-A}^{(d_r)}(c)} \right], \\
    & \qquad \qquad r=1,\ldots,I.
     \endaligned \tag 4.15 $$
Let $n \in \Lambda$ be a set of indices for which \thetag{4.9} holds.
By \thetag{4.15} for $r=1$ (any other $r$ gives the same result),
\thetag{3.27} and \thetag{3.28} (and its consequence, the $n$th root
asymptotics), one obtains by taking limits for $n \in \Lambda$ and
$n \to \infty$
$$   0 = h_\Lambda\left(\frac{\varphi(c)}2 \right) . $$
Let us now show that $\varphi(c)/2$ is a zero of $h_\Lambda$ of multiplicity
$I=I_j$. If the indices $d_1,\ldots,d_I$ were consecutive numbers, then
Leibniz' rule of derivation, via using \thetag{4.15}, would give us this rather
directly (see \cite{8} where $d_1=0, d_2=1,\ldots,d_I=I-1$). Here the case
is a bit more complicated. To see that $h_\Lambda'(\varphi(c)/2) = 0$ we
work out the following formula. Let $p,q$ be arbitrary functions
for which all operations we are about to carry out are admissible. This will be
true for our functions for $n \geq n_0$. Then
$$  \frac{q}{q'} \left( \frac{p}{q} \right)' = \left( \frac{p'}{q'} -
   \frac{p}{q} \right) , $$
and thus for any $\tau \in {\Bbb N}$
$$  \frac{q^{(\tau-1)}}{q^{(\tau)}}
   \left( \frac{p^{(\tau-1)}}{q^{(\tau-1)}} \right)'
     = \frac{p^{(\tau)}}{q^{(\tau)}} -
      \frac{p^{(\tau-1)}}{q^{(\tau-1)}} . $$
Because of the structure of the right hand side, by summing up we obtain
$$      \frac{p^{(\tau)}}{q^{(\tau)}} - \frac{p}{q} =
    \frac{q}{q'} \left( \frac{p}{q} \right)'
    + \frac{q'}{q''} \left( \frac{p'}{q'} \right)'
    + \cdots +
    \frac{q^{(\tau-1)}}{q^{(\tau)}} \left( \frac{p^{(\tau-1)}}{q^{(\tau-1)}}
   \right)' .  \tag 4.16 $$
Taking $p=q_{n-k}^{(d_1)}$, $q=q_{n-A}^{(d_1)}$ and $\tau=d_2-d_1$, from
\thetag{4.15} and \thetag{4.16} follows that
(multiply the difference for $r=1$ and $r=2$ by $n$)
$$  \multline
   0 = \sum_{k=0}^A \lambda_n^* \lambda_{n,k}
     \frac{\kappa_{n-A}}{\kappa_{n-k}}
     \left[ O \left( \frac{n^{N+1}}{q_{n-A}(c)q_{n-A}^{(d_1)}(c)} \right)
     \right. \\
    + \left. \sum_{\tau=0}^{d_2-d_1-1}
    \frac{n q_{n-A}^{(d_1+\tau)}(c)}{q_{n-A}^{(d_1+\tau+1)}(c)}
    \left( \frac{q_{n-k}^{(d_1+\tau)}}{q_{n-A}^{(d_1+\tau)}} \right)'(c)
    \right] . \endmultline $$
Now, from \thetag{3.27}--\thetag{3.29} and \thetag{4.9}, taking limits
for $n \in \Lambda$ and $n \to \infty$, we get
$$  0 = \sqrt{c^2-1} (d_2-d_1) (h_\Lambda \circ \varphi/2)'(c), $$
or, equivalently,
$$  0 = h_\Lambda'\left( \frac{\varphi(c)}2 \right) , $$
since $\varphi'(c) \neq 0$. This procedure can be continued until we use all
the $I$ equations in \thetag{4.15}. We show how this is done. Leibniz'
formula gives
$$  0 = - \left( \frac{p^{(\eta)}}{q^{(\eta)}} \right)^{(\ell)}
     + \sum_{\nu=0}^{\ell} {\ell \choose \nu}
     \frac{p^{(\eta+\ell-\nu)}}{q^{(\eta+\ell-\nu)}}
     \left( \frac{1}{q^{(\eta)}} \right)^{(\nu)}
     q^{(\eta+\ell-\nu)}. \tag 4.17 $$
Take $\eta$ from $d_1$ up to $d_{\ell+1}-\ell$. We obtain a total of
$d_{\ell+1}-\ell-d_1+1$
equations. In these equations, the quotients
$p^{(\eta+\ell-\nu)}/q^{(\eta+\ell-\nu)}$ run from
$p^{(d_1)}/q^{(d_1)}$ up to $p^{(d_{\ell+1})}/q^{(d_{\ell+1})}$.
Of these quotients, we assume $p^{(d_r)}/q^{(d_r)}$, $r=1,\ldots,\ell+1$,
known, and the rest of them, unknown. In this
way, we have a total of $d_{\ell+1}-d_1+1-(\ell+1) = d_{\ell+1}-d_1-\ell$
unknowns. We add the fictitious `unknown' $-1$, the coefficient of
$(p^{(\eta)}/q^{(\eta)})^{(\ell)}$.
Put the terms corresponding to the known quotients on one side of the
equations,
then we have a system of linear equations with the same number
of equations and unknowns.
Denote by $D_0$ the determinant of this linear system, and by $D_{-1}$
the determinant corresponding to the `unknown' $-1$, then by Cramer's rule
$$  D_0 = -D_{-1}.  \tag 4.18 $$
Expanding the
determinants $D_0$ and $D_{-1}$ along their first columns one sees that
formula \thetag{4.18} expresses a certain linear combination of
the expressions
$(p^{(\eta)}/q^{(\eta)})^{(\ell)}$, $\eta=1,\ldots,d_{\ell+1}-\ell$, as
a linear combination of
$p^{(d_r)}/q^{(d_r)}$, $r=1,\ldots,\ell+1$,
say
$$  \sum_{r=1}^{\ell+1} \beta_r \frac{p^{(d_r)}}{q^{(d_r)}}
   = \sum_{\eta=d_1}^{d_{\ell+1}-\ell}
     (-1)^\eta \alpha_\eta \left(
\frac{p^{(\eta)}}{q^{(\eta)}}\right)^{(\ell)}. \tag 4.19 $$
Note that the coefficients $\beta_r$ and $\alpha_\eta$,
only depend on the denominator $q$ and not on $p$ (see \thetag{4.17}). Each
$\beta_r$ is the sum of
products of $d_{\ell+1} - d_1 - \ell+1$ factors of the form
$$  \left( \frac{1}{q^{(\eta)}} \right)^{(\nu)} q^{(\eta+\ell-\nu)} . $$
On the other hand $\alpha_\eta$ is the determinant of the minor
corresponding to
$(p^{(\eta)}/q^{(\eta)})^{(\ell)}$ in $D_0$.
For each fixed $n \in  \Lambda$, take $q = q_{n-A}$ then $\beta_r =
\beta_r(n)$ and
$\alpha_\eta = \alpha_\eta(n)$. By Lemma 3 and \thetag{3.29}
$$ \multline
   \left( \frac{1}{q_{n-A}^{(\eta)}} \right)^{(\nu)}
     \frac{q_{n-A}^{(\eta+\ell-\nu)}}{n^\ell} =
     \left[ g_{n-A,\nu} \frac{n^\nu}{q_{n-A}^{(\eta)}}
       + o\left( \frac{n^\nu}{q_{n-A}^{(\eta)}} \right) \right]
       \frac{q_{n-A}^{(\eta+\ell-\nu)}}{n^\ell} \\
       \uto \left( \frac{-1}{\sqrt{z^2-1}} \right)^\nu
            \left( \frac{1}{\sqrt{z^2-1}} \right)^{\ell-\nu}
        = (-1)^\nu \left( \frac{1}{\sqrt{z^2-1}} \right)^\ell,
        \qquad K \subset {\Bbb C} \setminus S_\mu .
        \endmultline   \tag 4.20 $$
Therefore
$$   \beta_r = O \left(n^{\ell(d_{\ell+1}-d_1-\ell+1)} \right) .
    \tag 4.21 $$
Multiply equation \thetag{4.15} by $\beta_r(n)$, $r=1,\ldots,\ell+1$, and sum
these first $\ell+1$ formulas. From \thetag{4.19} and \thetag{4.21} we
obtain
$$  \multline
  0 = \sum_{k=0}^A \lambda_n^* \lambda_{n,k} \frac{\kappa_{n-A}}{\kappa_{n-k}}
    \left[ O \left( \frac{n^{N+\ell(d_{\ell+1}-d_1-\ell+1)}}{q_{n-A}(c)
     q_{n-A}^{(d_{\ell+1})}(c)} \right) \right. \\
     + \left.
     \sum_{\eta=d_1}^{d_{\ell+1}-\ell} (-1)^\eta \alpha_\eta(n)
     \left( \frac{q_{n-k}^{(\eta)}}{q_{n-A}^{(\eta)}} \right)^{(\ell)}(c)
     \right] .
   \endmultline \tag 4.22 $$
In deriving \thetag{4.22} we recall that the coefficients $\beta_r(n)$
and $\alpha_\eta(n)$ do not depend on the numerator $p$,
therefore \thetag{4.19} is also valid
with the same coefficients for  $p=q_{n-k}$, $k=0,\ldots,A$.
From \thetag{4.20}
we also know that $\alpha_\eta(n)$, for each $\eta$, is of the order
$n^{\ell(d_{\ell+1}-d_1-\ell)}$.
Therefore, dividing \thetag{4.22} by this quantity, taking limits for
$n \to \infty,\ n \in \Lambda$, we obtain from
\thetag{3.27}--\thetag{3.29}, \thetag{4.9} and \thetag{4.20}
$$ 0 = \lim_{n \to \infty,\ n \in \Lambda}
  \left[ \sum_{\eta=d_1}^{d_{\ell+1}-\ell}
   \frac{(-1)^\eta \alpha_\eta(n)}{n^{\ell(d_{\ell+1}-d_1-\ell)}} \right]
    h_\Lambda^{(\ell)} \left( \frac{\varphi(c)}2 \right) =
    \frac{M h_\Lambda^{(\ell)} \left( \frac{\varphi(c)}2 \right)}
    {(\sqrt{c^2-1})^{\ell(d_{\ell+1}-d_1-\ell)}} , $$
where $M \neq 0$ is the determinant of the matrix obtained by deleting columns
1, $d_2 - d_1+1,\ldots,d_{\ell+1} - d_1 + 1$ (counted from right to left) of the
matrix
$$ \matrix
        & k-1 \text{ columns} \vspace{3pt} \\
      d_{\ell+1}-d_1-\ell+1 \text{ rows} &
\pmatrix
  1 & 0 & \cdots & 0 & 0 & {\ell \choose 0} & - {\ell \choose 1} &
         \cdots & (-1)^\ell {\ell \choose \ell} \\
  1 & 0 & \cdots & 0 & {\ell \choose 0} & - {\ell \choose 1} & {\ell \choose 2}
          & \cdots & 0 \\
  \vdots & \vdots & \cdots &  & & & & \cdots & \vdots \\
  1 & {\ell \choose 0} & \cdots & & & & & \cdots & 0
  \endpmatrix
  \endmatrix  $$
(see \thetag{4.17}).
This matrix is the result of shifting the row
${\ell \choose 0}, -{\ell \choose 1}, \ldots, (-1)^\ell {\ell \choose \ell}$
to the left $k-1$ times and padding the remaining part with zeros, but with a
column of one's on the extreme left.
Since this determinant is nonzero it follows that
 $h_\Lambda(z)$ is divisible  by $(z-\varphi(c_j)/2)^{I_j}$ for each
$j=1,\ldots,m$.

Now we use the $N+1-I$ equations in \thetag{4.11} for those indices $i$ for
which $\L_i \equiv 0$.
 Before doing this let us look more closely at equations
\thetag{4.10} in order to derive one more property of the polynomials $p_i$. If
one writes down the equations \thetag{4.10}
one finds (again we drop the index $j$)
$$ \align
  0 &= \int S_n \frac{d\rho}{x-c} + N! \frac{s}{0!} \L_N \\
  0 &= \int S_n \frac{d\rho}{(x-c)^2} + N! \frac{s'}{1!} \L_N + (N-1)!
    \frac{s}{0!} \L_{N-1} \\
  0 &= \int S_n \frac{d\rho}{(x-c)^3} + N! \frac{s''}{2!} \L_N +
  (N-1)! \frac{s'}{1!} \L_{N-1} + (N-2)! \frac{s}{0!} \L_{N-2} \\
  0 &= \int S_n \frac{d\rho}{(x-c)^4} + N! \frac{s^{(3)}}{3!} \L_N + (N-1)!
 \frac{s''}{2!} \L_{N-1} + (N-2)! \frac{s'}{1!} \L_{N-2} + \\
   &\ \qquad (N-3)! \frac{s}{0} L_{N-3} \\
 \vdots &\  \qquad\qquad \vdots \\
 0 &= \int S_n \frac{d\rho}{(x-c)^{N+1}} + N! \frac{s^{(N)}}{N!} \L_N +
 (N-1)!  \frac{s^{(N-1)}}{(N-1)!} \L_{N-1} + \cdots  \\
  &\ \qquad + 1! \frac{s'}{1!} \L_1 + 0! \frac{s}{0!} \L_0 .
 \endalign $$
From this, finding $p_N$ is immediate. To obtain $p_{N-1}$ we take
$\beta_1'=1$ and $\beta_2'$ satisfying
$$  \beta_1' \frac{s'}{1!} + \beta_2' s = 0 $$
and from the first two equations obviously
$$  0 = \int S_n \left( \frac{\beta_1'}{(x-c)^2} + \frac{\beta_2'}{x-c}
    \right) \, d\rho + (N-1)! s \L_{N-1}, $$
from which  $p_{N-1}$ is determined. Taking $\beta_1'$ and
$\beta_2'$ as above, we then take $\beta_3'$ such that
$$  \beta_1' \frac{s''}{2!} + \beta_2' \frac{s'}{1!} + \beta_3' s = 0, $$
and from the first three equations one sees that
$$  0 = \int S_n \left( \frac{\beta_1'}{(x-c)^3} + \frac{\beta_2'}{(x-c)^2}
   + \frac{\beta_3'}{x-c} \right) \, d\rho + (N-2)! \frac{s}{0!} \L_{N-2}, $$
from which we can find $p_{N-2}$. Continuing this process we obtain all
$p_i$ (and $\beta_1',\ldots,\beta_{N+1}')$.

Let $i_1 < \cdots < i_{N+1-I}$ denote those indices $i$ for which $\L_i
 \equiv  0$. Note that
$$  \text{degree } p_{i_r} = d_r' = N+1-i_r \geq 1. $$
Therefore all  $p_{i_r}$ are of different degree and their
coefficients do not depend on $n$ and moreover their general form is such
that
$$  p_{i_r} \left( \frac{1}{x-c} \right) = \sum_{\delta=1}^{d_r'}
   \frac{\beta_\delta'}{(x-c)^{d_r+1-\delta}} .  \tag 4.23 $$
As before, we show how the proof works in the first
two steps, and then give the general outline on how to carry on.
If we multiply \thetag{4.11} for $i_r$ by $q_{n-A}^{(d_1-d_r)}(c)$ and
use orthogonality (in the second equality we use the same arguments
as those we used to obtain \thetag{4.13}, and then in the third equality we
eliminate other terms also equal to zero by orthogonality) we obtain
for $r=1,\ldots,N+1-I$
$$  \aligned
  0 &= q_{n-A}^{(d_1'-d_r')}(c) \int S_n(x) p_r \left( \frac{1}{x-c} \right)
  \, d\rho(x) \\
 &= \int S_n(x) \sum_{\eta=0}^{d_r'-1} \frac{q_{n-A}^{(d_1'-d_r'+\eta)}}{\eta!}
    (c-x)^\eta p_{i_r} \left( \frac{1}{x-c} \right) \, d\rho(x) \\
 &= \sum_{\delta=1}^{d_r'} \beta_\delta' \int S_n(x) \left[
   \sum_{\eta=0}^{d_r'-\delta} (-1)^\eta
   \frac{q_{n-A}^{(d_1'-d_r'+\eta)}}{\eta!}
   \frac{1}{(x-c)^{d_r'+1-\delta-\eta}} \right] \, d\rho(x) .
   \endaligned      \tag 4.24  $$
Now it is easy to prove that if $\Lambda \subset {\Bbb Z}_+$ is such that
\thetag{4.9} holds, then
$$   h_\Lambda \left( \frac{1}{2\varphi(c)} \right) = 0 . $$
In fact, note that in \thetag{4.24} the
highest order of derivative appearing is $d_1' - 1$ (for $\delta =1$ and
$\eta =d_r-1$). If we multiply any one of these equations (say for $r=1$) by
$\lambda_n^* \kappa_{n-A}/n^{d_1'-1}$,
using \thetag{4.6}, \thetag{3.33} and \thetag{2.13}, we obtain by taking the
limit for $n \to \infty,\ n \in \Lambda$, that
$$  \left( \frac{1}{\sqrt{c^2-1}} \right)^{d_1'-1} h_\Lambda \left(
    \frac{1}{2\varphi(c)} \right) = 0 . $$
In order to prove that
$h_\Lambda' \left( \frac{1}{2\varphi(c)} \right) = 0$
we need to eliminate
the term with derivative $d_1'-1$. In order to do this we combine the first two
equations. Note that  $d_1' > d_2' \geq 1$ (actually $d_2 \geq  2$
because $\L_N \not\equiv  0$). Since
$$  \vmatrix \frac{1}{(d_1'-1)!} & \frac{1}{(d_1'-2)!} \\
             \frac{1}{(d_2'-1)!} & \frac{1}{(d_2'-2)!}
     \endvmatrix \neq 0, $$
there exist $\alpha_1$ and $\alpha_2$ such that
$$   \alpha_1 \frac{(-1)^{d_1'-1}}{(d_1'-1)!}
    + \alpha_2 \frac{(-1)^{d_2'-1}}{(d_2'-1)!} = 0 , $$
and
$$   \alpha_1 \frac{(-1)^{d_1'-2}}{(d_1'-2)!}
    + \alpha_2 \frac{(-1)^{d_2'-2}}{(d_2'-2)!} = 1 . $$
Multiply the first equation in \thetag{4.24} $(r=1)$ by
$\alpha_1 \lambda_n^* \kappa_{n-A}/n^{d_1'-2}$,
 the second one $(r=2)$ by
$\alpha_2 \lambda_n^* \kappa_{n-A}/n^{d_1'-2}$,
 sum them up and take the limit for $n \to \infty,\ n \in  \Lambda$. Again
using \thetag{4.6},
\thetag{3.33} and \thetag{2.13}, we obtain (even if $\beta_2 \neq  0$
because $h_\Lambda \left( \frac{1}{2\varphi(c)} \right) = 0$)
$$  0 = \left[ \left( \frac{1}{\sqrt{z^2-1}} \right)^{d_1'-1} h_\Lambda
   \left( \frac{1}{2\varphi} \right) \right]'(c)
    + \beta_2 \left( \frac{1}{\sqrt{c^2-1}} \right)^{d_1'-1}
    h_\Lambda \left( \frac{1}{2\varphi(c)} \right) . $$
Using  $h_\Lambda \left( \frac{1}{2\varphi(c)} \right) = 0$
this formula implies
$$  h_\Lambda' \left( \frac{1}{2\varphi(c)} \right) = 0 . $$
The proof is completed by induction using the scheme employed in passing to the
first derivative. Of course the situation becomes more complicated since each
time there appear
more terms with the same order of derivative. Fortunately, this causes no real
problem because the polynomials $p_{i_r}$ have
correspondingly equal coefficients starting from the highest
degree down (see \thetag{4.23}).
To be more precise, in proving that
$$  h_\Lambda^{(\tau-1)} \left( \frac{1}{2\varphi(c)} \right) = 0
  \qquad \tau \leq N+1-I, $$
we start with the induction hypothesis
that
$$  h_\Lambda^{(\tau')} \left( \frac{1}{2\varphi(c)} \right) = 0
  \qquad \tau'=0,\ldots,\tau-2. $$
Noting that
$d_1'>d_2'> \cdots >d_\tau' \geq 1$
(actually $\geq 2$) we will show that
$$  \vmatrix
  \frac{1}{(d_1'-1)!} & \frac{1}{(d_1'-2)!} & \cdots & \frac{1}{(d_1'-\tau)!} \\
  \frac{1}{(d_2'-1)!} & \frac{1}{(d_2'-2)!} & \cdots & \frac{1}{(d_2'-\tau)!} \\
   \vdots & \vdots & \cdots & \vdots \\
  \frac{1}{(d_\tau'-1)!} & \frac{1}{(d_\tau'-2)!} & \cdots &
\frac{1}{(d_\tau'-\tau)!}
   \endvmatrix \neq 0 . $$
Indeed, this statement is obviously equivalent to
$$  \vmatrix
     1 & d_1'-1 & (d_1'-1)(d_1'-2) & \cdots \\
     1 & d_2'-1 & (d_2'-1)(d_2'-2) & \cdots \\
     \vdots & \vdots & \vdots & \cdots \\
     1 & d_\tau'-1 & (d_\tau'-1)(d_\tau'-2) & \cdots
     \endvmatrix \neq 0 . $$
This determinant can be transformed into a Vandermonde determinant by
elementary column operations. Therefore, there exist $\alpha_1, \ldots,
\alpha_\tau$ such that
$$  \sum_{r=1}^\tau  \alpha_r \frac{(-1)^{d_r'-k}}{(d_r'-k)!} = 0
   \qquad k=1,\ldots,\tau-1, $$
and
$$  \sum_{r=1}^\tau  \alpha_r \frac{(-1)^{d_r'-\tau}}{(d_r'-\tau)!} = 1. $$
Take the first $\tau$ equations in \thetag{4.24}, multiply the $r$th equation
by  $\alpha_r \lambda_n^* \kappa_{n-A}/n^{d_1'-\tau}$,
sum them up and observe that in the sum all
derivatives of order $d_1'-1,\ldots,d_1'-\tau+1$  cancel out. Using
\thetag{4.6},
\thetag{3.33} and \thetag{2.13} and taking the limit for $n \to \infty,\ n \in
\Lambda$ then gives
$$   0 = \sum_{r=1}^\tau \beta_r' \left[ \left( \frac{1}{\sqrt{z^2-1}}
   \right)^{d_1'-\tau+1} h_\Lambda \left( \frac{1}{2\varphi} \right)
   \right]^{(\tau-r)}(c), $$
and using the induction hypothesis we immediately obtain that
$$   h_\Lambda^{(\tau-1)} \left( \frac{1}{2\varphi(c)} \right) = 0 . $$
Hence $h_\Lambda(z)$ is also divisible by
$\left( z- \frac{1}{2\varphi(c)} \right)^{N+1-I}$ for each $j=1,\ldots,m$.
 Therefore \thetag{4.7} holds and thus
$$  \frac{S_n(z)}{Q_{n-A}(z)} \uto
  \prod_{j=1}^m \left( \frac{\varphi(z)}2 - \frac{\varphi(c_j)}2
    \right)^{I_j} \left( \frac{\varphi(z)}2 - \frac{1}{2\varphi(c_j)}
     \right)^{N_j+1-I_j}, \qquad K \subset \overline{\Bbb C} \setminus S_\mu, $$
or, equivalently on account of \thetag{3.3} (see Lemma 4)
$$  \frac{S_n(z)}{Q_n(z)} \uto \prod_{j=1}^m
  \left( 1 - \frac{\varphi(c_j)}{\varphi(z)} \right)^{I_j}
    \left( 1 - \frac{1}{\varphi(z)\varphi(c_j)} \right)^{N_j+1-I_j},
    \qquad K \subset \overline{\Bbb C} \setminus S_\mu . \tag 4.25 $$
Finally, from Lemma 4 we also have that
$$  \frac{Q_n(z)}{L_n(z)} \uto \prod_{j=1}^m \left(
    \frac{\varphi(z)-\varphi(c_j)}{2(z-c_j)} \right)^{N_j+1}, \qquad
    K \subset \overline{\Bbb C} \setminus S_\mu .  \tag 4.26 $$
It is easy to verify that
$$  \left( \frac{\varphi(z)-\varphi(c)}{2(z-c)} \right)
   \left( 1 - \frac{1}{\varphi(z)\varphi(c)} \right) = 1,
   \qquad z,c \in {\Bbb C} \setminus [-1,1], $$
hence \thetag{4.25} and  \thetag{4.26} give
$$  \frac{S_n(z)}{L_n(z)} \uto
   \prod_{j=1}^m \left( 1 - \frac{\varphi(c_j)}{\varphi(z)} \right)^{I_j}
     \left( 1 - \frac{1}{\varphi(z)\varphi(c_j)} \right)^{-I_j},
     \qquad K \subset \overline{\Bbb C} \setminus S_\mu , $$
which is the same as \thetag{4.4} (for $\nu  = 0$).
For arbitrary $\nu$ as in Theorem 1, \thetag{4.4} follows by induction. The
statements concerning the degree of $S_n$ and the asymptotic behavior of
its zeros follow from \thetag{4.4} and Hurwitz' theorem. \qed
\enddemo

Formula \thetag{4.4} and the existence of the limits of the coefficients
$\lambda_{n,k}$ allow us to obtain, using arguments similar to those given
above, other asymptotic formulas for the polynomials $S_n$. We collect some of
them in the following corollary, the proof of which we leave to the reader.

\proclaim{Corollary 2}
Under the hypothesis of Theorem 4, for all sufficiently
large $n$,
$\langle S_n, S_n \rangle> \neq  0$. Let $\tau_n$ be the leading
coefficient of
$l_n$, the $n$th orthonormal polynomial with respect to $\mu$. Then $\gamma_n
= \langle S_n,S_n \rangle^{-1/2}$, $n\geq n_0$, may be taken so that
$$   \lim_{n \to \infty} \frac{\gamma_n}{\tau_n} =
   \prod_{j=1}^m \frac{1}{(\varphi(c_j))^{I_j}}, $$
and in particular
$$   \lim_{n \to \infty} \frac{\gamma_{n+1}}{\gamma_n} = 2 . $$
Denote by $s_n = \gamma_n S_n$, $n \geq n_0$,
the $n$th orthonormal polynomial with respect to $\langle \cdot,\cdot
\rangle$.  Then for all fixed $\nu \in {\Bbb Z}_+$
$$  \frac{s_n^{(\nu)}(z)}{\ell_n^{(\nu)}(z)} \uto
   \prod_{j=1}^m \left( \frac{(\varphi(z)-\varphi(c_j))^2}{2\varphi(z)
     \varphi(c_j) (z-c_j)} \right)^{I_j}, \qquad
     K \subset \overline{\Bbb C} \setminus S_\mu , $$
$$  \frac{s_{n+1}^{(\nu)}(z)}{s_n^{(\nu)}(z)} \uto
    \varphi(z), \qquad  K \subset \overline{\Bbb C} \setminus
    (S_\mu \cup \{ c_1,\ldots,c_m\} ), $$
$$  \frac1n \frac{s_n^{(\nu+1)}(z)}{s_n^{(\nu)}(z)} \uto
  \frac{1}{\sqrt{z^2-1}}, \qquad K \subset \overline{\Bbb C} \setminus
    (S_\mu \cup \{ c_1,\ldots,c_m\} ) . $$
\endproclaim

Theorem 4 and the asymptotic properties of orthogonal polynomials with respect
to measures in $\MC(0,1)$ allow us to obtain an extension of
Gonchar's result in \cite{8} (see also \cite{14}). Using the approach given in
\cite{10} one easily obtains

\proclaim{Corollary 3}
Let $f$ be as in \thetag{1.5} where $\mu \in  \MC(0,1)$. Let
$\pi_n$
be the $[n-1,n]$ Pad\'e approximant for $f$, $Q_n$ its
denominator with leading coefficient equal to 1 and $L_n$ the $n$th monic
orthogonal polynomial with respect to $\mu$. Then
$$   \frac{Q_n(z)}{L_n(z)} \uto
   \prod_{j=1}^m \left( \frac{(\varphi(z)-\varphi(c_j))^2}{2\varphi(z)
      (z-c_j)} \right)^{N_j+1}, \qquad
    K \subset \overline{\Bbb C} \setminus S_\mu . $$
For all sufficiently large $n$ the degree of $Q_n$ is $n$, and each $c_j$
attracts
exactly $N_j+1$ zeros of $Q_n$, while the rest of the zeros of $Q_n$ accumulate
on $S_\mu$. The rate of convergence of $\pi_n$ to $f$ is given by
$$  \frac{f(z) - \pi_{n+1}(z)}{f(z)-\pi_n(z)} \uto \frac{1}{\varphi^2(z)},
   \qquad  K \subset \overline{\Bbb C} \setminus
    (S_\mu \cup \{ c_1,\ldots,c_m\} )  . $$
In particular
$$   \lim_{n \to \infty} \| f-\pi_n \|_K^{\frac1{2n}} =
   \frac{1}{\|\varphi \|_K},    \qquad
   K \subset {\Bbb C} \setminus
    (S_\mu \cup \{ c_1,\ldots,c_m\} )   . $$
\endproclaim

As a final remark we wish to say that for regular generalized Sobolev inner
products it is not difficult to prove that the zeros attracted by the $c_j$'s
converge geometrically when $\mu \in \MC(0,1)$.
 Special cases show that in general this is not so if
the inner product is not regular. Nevertheless, it is plausible that formula
\thetag{4.4} remains valid in the nonregular case if we substitute $I_j$ by
the rank of $\Gamma_j^*$.

\head Acknowledgments \endhead
The authors wish to thank the referees for their constructive comments
and guidance. Their reports have led to a significantly improved manuscript.


\Refs
\ref \no 1
\by L.V. Ahlfors
\book Complex Analysis
\bookinfo (3rd edition)
\publ McGraw Hill \publaddr New York \yr 1979
\endref
\ref \no 2
\by  M. Alfaro, F. Marcell\'an and M.L. Rezola
\paper Orthogonal polynomials on Sobolev spaces: old and new directions
\jour J. Comput. Appl. Math. \vol 48 \yr 1993 \pages 113--132
\endref
\ref \no 3
\by O. Blumenthal
\book \"Uber die Entwicklung einer willk\"urlichen Funktion nach den
Nennern des Kettenbruches f\"ur $\int_{-\infty}^0 \frac{\varphi(\xi) \, d\xi}
{z-\xi}$
\publ Inaugural Dissertation \publaddr G\"ottingen \yr 1898
\endref
\ref \no 4
\by  T.S. Chihara
\book An Introduction to Orthogonal Polynomials
\publ Gordon and Breach \publaddr New York \yr 1978
\endref
\ref \no 5
\by W.D. Evans, L.L. Littlejohn, F. Marcell\'an, C. Markett and A. Ronveaux
\paper On recurrence relations for Sobolev orthogonal polynomials
\jour SIAM J. Math. Anal. \toappear
\endref
\ref \no 6
\by T.W. Gamelin
\book Uniform Algebras
\publ Prentice Hall \publaddr Englewood Cliffs, NJ \yr 1969
\endref
\ref \no 7
\by  J. Gilewicz
\paper Location of the zeros of polynomials satisfying three-term recurrence
relations with complex coefficients, I. General case
\jour J.  Approx. Theory \vol 43 \yr 1985 \pages 1--14
\endref
\ref \no 8
\by A.A. Gonchar
\paper On the convergence of Pad\'e approximants for
     some classes of meromorphic functions
\jour Mat. Sb. \vol 97 {\rm (139)} \yr 1975 \pages 607--629
\transl \jour Math. USSR Sb. \vol 26 \yr 1975 \pages 555--575
\endref
\ref \no 9
\by  P.P. Korovkin
\paper The capacity of a set and polynomials which minimize an integral
\jour Kaliningrad Ped. Inst. Uchenie Zap. \vol V \yr 1958 \pages 34--52
\endref
\ref \no 10
\by G. L\'opez
\paper Convergence of Pad\'e approximants of Stieltjes type meromorphic
functions and comparative asymptotics for orthogonal polynomials
\jour Mat. Sb. \vol 136 {\rm (178)} \yr 1988 \pages 46--66
\transl \jour Math. USSR Sb. \vol 64 \yr 1989 \pages 207--227
\endref
\ref \no 11
\by  G. L\'opez
\paper Relative asymptotics of orthogonal polynomials on the real axis
\jour Mat. Sb. \vol 137 {\rm (179)} \yr 1988 \pages 500--525
\transl \jour Math. USSR Sb. \vol 65 \yr 1990 \pages 505--527
\endref
\ref \no 12
\by  G. L\'opez
\paper Asymptotics of polynomials orthogonal with respect to varying measures
\jour Constr. Approx. \vol 5 \yr 1989 \pages 199--219
\endref
\ref \no 13
\by G. L\'opez, D. Barrios and E. Torrano
\paper Zero distribution and asymptotics of polynomials satisfying
 three-term recurrence relations with complex coefficients
\jour Mat. Sb. \toappear
\endref
\ref \no 14
\by A. Magnus
\paper Toeplitz matrix techniques and convergence of complex weight
 Pad\'e approximants
\jour J.  Comput.  Appl. Math. \vol 19 \yr 1987 \pages 23--38
\endref
\ref \no 15
\by F. Marcell\'an and W. Van Assche
\paper Relative asymptotics for orthogonal polynomials with a Sobolev inner
product
\jour J. Approx. Theory \vol 72 \yr 1993 \pages 193--209
\endref
\ref \no 16
\by A. M\'at\'e, P. Nevai and V. Totik
\paper Asymptotics for the ratio of leading coefficients of orthonormal
polynomials on the unit circle
\jour Constr. Approx. \vol 1 \yr 1985 \pages 63--69
\endref
\ref \no 17
\by A. M\'at\'e, P. Nevai and V. Totik
\paper Extensions of Szeg\H{o}'s theory of orthogonal polynomials, II; III
\jour Constr. Approx. \vol 3 \yr 1987 \pages 51--72; 73--96
\endref
\ref \no 18
\by A. M\'at\'e, P. Nevai and V. Totik
\paper Strong and weak convergence of orthogonal polynomials
\jour Amer. J. Math. \vol 109 \yr 1987 \pages 239--281
\endref
\ref \no 19
\by  P. Nevai
\book Orthogonal Polynomials
\bookinfo Memoirs Amer. Math. Soc. \vol 213
\publ Amer. Math. Soc. \publaddr Provindence, RI \yr 1979
\endref
\ref \no 20
\by E.A. Rakhmanov
\paper On the asymptotics of the ratio of orthogonal polynomials
\jour Mat. Sb. \vol 103 {\rm (145)} \yr 1977 \pages 271--291
\transl \jour Math. USSR Sb. \vol 32 \yr 1977 \pages 199--213
\endref
\ref \no 21
\by  E.A. Rakhmanov
\paper On the asymptotics of the ratio of orthogonal polynomials II
\jour Mat. Sb. \vol 118 {\rm (160)} \yr 1982 \pages 104--117
\transl \jour Math. USSR Sb. \vol \yr 1983 \vol 46 \pages 105--117
\endref
\ref \no 22
\by  E.A. Rakhmanov
\paper On the asymptotics of polynomials orthogonal on the circle with
weights not satisfying Szeg\H{o}'s condition
\jour Mat. Sb. \vol 130 {\rm (172)} \yr 1986 \pages 151--169
\transl \jour Math. USSR Sb. \vol 58 \yr 1987 \pages 149--167
\endref
\ref \no 23
\by W. Rudin
\book Real and Complex Analysis
\publ McGraw-Hill
\publaddr New York \yr 1966
\endref
\ref \no 24
\by J. Shohat
\book Th\'eorie G\'en\'erale des Polynomes Orthogonaux de Tchebichef
\bookinfo M\'emorial des Sciences Math\'ematiques \vol 66
\publ Gauthier-Villars \publaddr Paris \yr 1934
\endref
\ref \no 25
\by H. Stahl and V. Totik
\book General Orthogonal Polynomials
\publ Cambridge University Press \yr 1992
\endref
\ref \no 26
\by  J.L. Ullman
\paper On the regular behavior of orthogonal polynomials
\jour Proc. London Math. Soc. \vol 24 \yr 1972 \pages 119--148
\endref
\ref \no 27
\by E.B. Van Vleck
\paper On the convergence of algebraic continued fractions whose coefficients
have limiting values
\jour Trans. Amer. Math. Soc. \vol 5 \yr 1904 \pages 253--262
\endref
\ref \no 28
\by H. Widom
\paper Polynomials associated with measures in the complex plane
\jour J. Math. Mech. \vol 16 \yr 1967 \pages 997--1013
\endref
\endRefs
\enddocument